\RequirePackage{fix-cm}
\documentclass{svjour3_arxiv}                    
\smartqed  
\usepackage{amssymb}                   
\usepackage{amsbsy}                     
\usepackage{amsmath}
\usepackage{cases}
\usepackage{graphicx} 
\usepackage{import}						
\usepackage{xcolor}						
\usepackage{subcaption}     				
\usepackage{adjustbox}
\usepackage{booktabs}
\usepackage{multirow}
\usepackage{tikz}
\usepackage{pgfplots}
\usepgfplotslibrary{external}
\tikzsetexternalprefix{img/tikz/ext/ext_} 
\newlength\figureheight
\newlength\figurewidth 
\usepackage{algorithm}
\usepackage[noend]{algpseudocode}
\usepackage{listings}
\usepackage[binary-units=true, detect-weight=true, detect-inline-weight=math]{siunitx}
\usepackage[hidelinks]{hyperref}
\usepackage{url}
\usepackage{xspace}
\usepackage{setspace}
\makeatletter
\DeclareRobustCommand\ttfamilyLst
        {\not@math@alphabet\ttfamily\mathtt
         \fontfamily\ttdefault\scriptsize\selectfont}
\makeatother
\newcommand{\vm}[1]{\boldsymbol{#1}}
\newcommand{\mb}[1]{\boldsymbol{#1}}
\newcommand{\ind}[1]{^{#1}} % Index oben

\newcommand{\dd}{{\rm d}}
\newcommand{\TT}{^\mathsf{T}}
\newcommand{\eye}[1]{\vm{I}_{\!#1}}

\newcommand{\diag}{{\rm diag}}

\newcommand{\zweinorm}[1]{\ensuremath \left\| #1 \right\|_2}
\newcommand{\opt}[1]{\ensuremath{{#1}^*}}

\newcommand{\Tab}[1]{\mbox{Table~#1}\xspace}
\newcommand{\Fig}[1]{\mbox{Figure~#1}\xspace}

\newcommand{\Sec}[1]{\mbox{Section~#1}\xspace}

\newcommand{\Matlab}{{\rm \textsc Matlab}\xspace}
\newcommand{\Simulink}{{\rm \textsc Simulink}\xspace}

\newcommand{\Dspace}{{\rm \textsc dSpace}\xspace}
\newcommand{\GRAMPC}{{\rm \textsc GRAMPC}\xspace}

\newcommand{\ACADO}{{\rm \textsc ACADO}\xspace}
\newcommand{\VIATOC}{{\rm \textsc VIATOC}\xspace}
\newcommand{\Renesas}{{\rm \textsc Renesas}\xspace}
\RequirePackage{color}
\definecolor{myred}{rgb}{1,0,0}
\definecolor{myblue}{rgb}{0,0,1}
\definecolor{mygreen}{rgb}{0,1,0}

\newcommand{\lag}{\mu}
\newcommand{\pen}{c}
\newcommand{\penmatrix}{C}

\newcommand{\costIII}[1]{\ensuremath{\bar{ #1}}}

\newcommand{\myalg}[2]{\hrule\vspace{0.5mm}\textbf{Algorithm~#1} #2\newline\vspace{-2.5mm}\hrule}

 \journalname{Optimization and Engineering}
\begin{document}
\title{A software framework for embedded
	nonlinear model predictive 
	control using a gradient-based augmented Lagrangian approach 
	(GRAMPC)}
\titlerunning{A software framework for embedded nonlinear MPC (GRAMPC)}
\author{Tobias Englert          \and
        Andreas V\"olz			\and
        Felix Mesmer			\and
        S\"onke Rhein			\and
        Knut Graichen
}
\institute{Institute of Measurement, Control, and Microtechnology, Ulm University, Albert-Einstein-Allee 41, 89081 Ulm, Germany. \email{$<$firstname$>$.$<$lastname$>$@uni-ulm.de}
}
\date{}
\maketitle
\vspace{-12mm}
\begin{abstract}
	A nonlinear MPC framework is presented that is
	suitable for dynamical systems with sampling times in the
	(sub)millisecond range and that allows for an efficient implementation on
	embedded hardware. The algorithm is based on an augmented Lagrangian
	formulation with a tailored gradient method for the inner
	minimization problem. The algorithm is implemented in the software
	framework GRAMPC and is a fundamental revision of an earlier
	version. Detailed performance results are presented for a test
        set of benchmark problems and in comparison to other nonlinear
        MPC packages. In addition, runtime results and memory
        requirements for GRAMPC on ECU level demonstrate its
        applicability on embedded hardware.
\end{abstract}
\section{Introduction}\label{sec:intro}

Model predictive control (MPC) is one of the most
popular advanced control methods due to its ability to handle linear
and nonlinear systems with constraints and multiple inputs. 
The well-known challenge of MPC is the numerical effort that
is required to solve the underlying optimal control problem (OCP) online.
In the recent past, however, the methodological as well as
algorithmic development of MPC for linear and nonlinear systems has
matured to a point that MPC nowadays can be applied to highly
dynamical problems in real-time.
Most approaches of real-time MPC either rely on suboptimal solution
strategies~\cite{Scokaert:TAC:1999,Diehl:IEE:2004,Graichen:TAC:2010}
and/or use tailored optimization algorithms to 
optimize the computational efficiency.

Particularly for linear systems, the MPC problem can be reduced to a quadratic 
problem, for which the optimal control over the admissible polyhedral set 
can be precomputed. This results in an explicit MPC strategy with
minimal computational effort for the online
implementation~\cite{Bemporad:TAC:2002,Bemporad:AUT:2002}, 
though this approach is typically limited to 
a small number of state and control variables.
An alternative to explicit MPC is the online active set
method~\cite{Ferreau:IJRNC:2008} that takes advantage of the receding
horizon property of MPC in the sense that typically only a small
number of constraints becomes active or inactive from one MPC 
step to the next. 

In contrast to active set strategies, interior point methods relax the
complementary conditions for the constraints and therefore solve a
relaxed set of optimality conditions in the interior of the admissable
constraint set. The MPC software packages
FORCES (PRO)~\cite{Domahidi:CDC:2012} and
fast\_mpc~\cite{Wang:CST:2010} 
employ interior point methods for linear MPC problems.
An alternative to active set and interior point methods are
accelerated gradient methods~\cite{Richter:TAC:2012} 
that originally go back to Nesterov's
fast gradient method for convex problems~\cite{Nesterov:2003}. A
corresponding software package for 
linear MPC is FiOrdOs~\cite{Jones:NPCC:2012}.

For nonlinear systems, one of the first real-time
MPC algorithms was the continuation/GMRES
method~\cite{Ohtsuka:AUT:2004} that solves the optimality conditions
of the underlying optimal control problem based on a continuation strategy. 
The well-known ACADO Toolkit~\cite{Houska2011} uses the
above-mentioned active set strategy in combination with a real-time
iteration scheme to efficiently solve nonlinear MPC problems. 
Another recently presented MPC toolkit is VIATOC~\cite{Kalmari2015}
that employs a  projected gradient method to solve the
time-discretized, linearized MPC problem. %optimization problem of the MPC. 

Besides the real-time aspect of MPC, a current focus of research is on
embedded MPC, i.e.\ the neat integration of MPC on embedded hardware
with limited ressources. This might be field 
programmable gate arrays (FPGA)
\cite{Ling:IFACWC:2008,Kaepernick:UKACC:2014,Hartley:OCAM:2015}, 
programmable logic controllers (PLC) in standard
automation systems~\cite{Kufoalor:MED:2014,Kaepernick:IFACWC:2014}
or electronic control units (ECU) in automotive
applications~\cite{Mesmer:CST:2018}.  
For embedded MPC, several challenges arise in addition to
the real-time demand. For instance, numerical robustness even
at low computational accuracy and tolerance against infeasibility are
important aspects as well as the ability to provide fast, possibly
suboptimal iterates with minimal computational demand. 
In this regard, it is also desirable to satisfy the system dynamics in
the single MPC iterations in order to maintain dynamically consistent
iterates. 
Low code complexity for portability and a small memory footprint are
further important aspects 
to allow for an efficient implementation of
MPC on embedded hardware.
The latter aspects, in particular, require the usage of streamlined,
self-contained code rather than highly complex MPC algorithms.
To meet these challenges, gradient-based algorithms become popular
choices due to their general simplicity and low computational
complexity~\cite{Giselsson:IFACWC:2014,Kouzoupis:EEC:2015,Necoara:SCL:2015}.

Following this motivation, this paper presents a software framework
for nonlinear MPC that can be efficiently used for embedded control of
nonlinear and highly dynamical systems with sampling times in the
(sub)millisecond range. 
The presented framework is a fundamental revision
of the MPC toolbox GRAMPC \cite{Kaepernick2014} 
(Gradient-Based MPC -- [gr\ae mp$'$si:]) 
that was originally developed for nonlinear systems with pure input
constraints. 
The revised algorithm of GRAMPC presented in this paper 
allows to account for general nonlinear 
equality and inequality constraints, as well as terminal constraints. Beside
``classical'' MPC, the toolbox can be applied to MPC on shrinking
horizon, general optimal control problems, moving horizon estimation
and parameter optimization problems including free end time problems. 
The new algorithm is based on an augmented Lagrangian formulation
in connection with a real-time gradient
method and tailored line search and multiplier update strategies that
are optimized for a time and memory efficient implementation on embedded
hardware. The performance and effectiveness of augmented Lagrangian
methods for embedded nonlinear MPC was recently demonstrated
for various application examples on rapid prototyping and 
ECU hardware level~\cite{Harder:AIM:2017,Mesmer:CST:2018,Englert:CEP:2018}.
Beside the presentation of the augmented Lagrangian algorithm and the
general usage of GRAMPC, the paper compares its performance to the
nonlinear MPC toolkits ACADO and VIATOC for different benchmark
problems. Moreover, runtime results are presented for GRAMPC on dSPACE
and ECU level including its memory footprint to demonstrate its
applicability on embedded hardware.

The paper is organized as follows. Section~\ref{sec:ProbForm} presents the general
problem formulation and exemplarily illustrates its application to
model predictive control and moving horizon estimation. 
Section~\ref{sec:OptAlg} describes the augmented Lagrangian
framework in combination with a gradient method for the inner
minimization problem. Section~\ref{sec:StrUsage} gives an overview on
the structure and usage of GRAMPC. 
Section~\ref{sec:performance} evaluates the performance
of GRAMPC for different benchmark problems and in comparison to ACADO
and VIATOC, before Section~\ref{sec:conclusions} closes the paper.

Some norms are used inside the paper, in particular in Section~\ref{sec:OptAlg}.
The Euclidean norm of a vector $\vm x\in \mathbb{R}^n$ is denoted by 
$\|\vm x\|_2$, the weighted quadratic norm by 
$\|\vm x\|_{\vm Q}=(\vm x\TT\vm Q\vm x)^{1/2}$ for some positive
definite matrix $\vm Q$, and the scalar product of two vectors $\vm x,\vm y\in\mathbb{R}^n$
is defined as $\langle\vm x,\vm y\rangle=\vm x\TT \vm y$. For a time function $\vm x(t)$,
$t\in[0,T]$ with $T<\infty$, the vector-valued $L^2$-norm is defined
by $\|\vm x\|_{L_2}=\big(\sum_{i=1}^n \| x_i \|_{L_2}\big)^{1/2}$
with $\| x_i \|_{L_2}=\big(\int_0^T\!x_i^2(t)\,\dd t\big)^{1/2}$.
The supremum-norm is defined componentwise in the sense of
$\|\vm x\|_{L_\infty}=\big[\|x_1\|_{L_\infty} \ldots \|x_n\|_{L_\infty} \big]\TT$
with $\|x_i\|_{L_\infty}=\sup_{t\in[0,T]} |x_i(t)|$.
The inner product is denoted by 
$\langle\vm x,\vm y\rangle=\int_0^T \!\vm x\TT(t)\vm y(t)\,\dd t$
using the same (overloaded) $\langle\,\rangle$-notation as in the vector
case.
Moreover, function arguments (such as time $t$) might be omitted
in the text for the sake of enhancing readability.

\section{Problem formulation}
\label{sec:ProbForm}
This section describes the class of optimal control problems that can be
solved by \GRAMPC.
The framework is especially suitable for model predictive control and
moving horizon estimation, as the numerical solution method is
tailored to embedded applications. Nevertheless, \GRAMPC can be used to solve
general optimal control problems or parameter optimization problems as
well. 
\subsection{Optimal control problem}
\GRAMPC solves nonlinear constrained optimal control problems with
fixed or free end time and potentially unknown parameters. 
Consequently, the most generic problem formulation that can be
adressed by \GRAMPC is given by  
\begin{subequations}\label{eq:OCP_orig}
\begin{alignat}{3}
\label{eq:OCP_cost} 
& \!\!\!\min_{\vm u, \vm p, T} \quad
&& J(\vm u, \vm p, T;\vm x_0) = V(\vm x(T), \vm p, T)
+ \int_0^T l(\vm x(t), \vm u(t), \vm p, t) \, \dd t 
\\ \label{eq:OCP_dynamics}
& \textrm{s.t.} &&
\vm M \vm{\dot x}(t) = \vm f(\vm x(t), \vm u(t), \vm p, t) \,,
\quad \vm x(0) = \vm x_0 
\\ \label{eq:OCP_eqconstr}
&&&  \vm g(\vm x(t), \vm u(t), \vm p, t) = \vm 0 \,, \quad
\vm g_T(\vm x(T), \vm p, T) = \vm 0
\\ \label{eq:OCP_ieqconstr}
&&&  \vm h(\vm x(t), \vm u(t), \vm p, t) \le \vm 0 \,, \quad
\vm h_T(\vm x(T), \vm p, T) \le \vm 0
\\
&&& \vm u(t) \in \left[\vm u_{\min}, \vm u_{\max}\right] 
&&
\label{eq:OCP_boxedUconst}
\\
&&& \vm p \in \left[\vm p_{\min}, \vm p_{\max}\right] \,,\quad 
T \in \left[T_{\min}, T_{\max}\right]
\label{eq:OCP_boxedpTconst}
\end{alignat}
\end{subequations}
with state $\vm x\in \mathbb{R}^{N_{\vm x}}$, control $\vm u \in
\mathbb{R}^{N_{\vm u}} $, parameters $\vm p\in \mathbb{R}^{N_{\vm p}}$ and 
end time $T\in \mathbb{R}$. 
The cost to be minimized~\eqref{eq:OCP_cost} consists of the terminal and
integral cost functions 
$V:\mathbb{R}^{N_{\vm x}}\times\mathbb{R}^{N_{\vm p}}\times\mathbb{R} \rightarrow
\mathbb{R}$ and
$l:\mathbb{R}^{N_{\vm x}}\times\mathbb{R}^{N_{\vm u}}\times\mathbb{R}^{N_{\vm p}}\times\mathbb{R}
\rightarrow \mathbb{R}$, respectively.  
The dynamics~\eqref{eq:OCP_dynamics} are given in semi-implicit form
with the (constant) mass matrix $\vm M$, the nonlinear system function
$\vm
f:\mathbb{R}^{N_{\vm x}}\times\mathbb{R}^{N_{\vm u}}\times\mathbb{R}^{N_{\vm p}}\times\mathbb{R}
\rightarrow \mathbb{R}^{N_{x}}$, and the initial state $\vm x_0$.  
The system class~\eqref{eq:OCP_dynamics} includes standard ordinary
differential equations for $\vm M=\vm I$ as well as
(index-1) differential-algebraic equations with singular mass matrix
$\vm M$.
In addition, \eqref{eq:OCP_eqconstr} and \eqref{eq:OCP_ieqconstr}
account for equality and inequality constraints 
$\vm
g:\mathbb{R}^{N_{\vm x}}\times\mathbb{R}^{N_{\vm u}}\times\mathbb{R}^{N_{\vm p}}\times\mathbb{R}
\rightarrow \mathbb{R}^{N_{\vm g}}$ and 
$\vm h:\mathbb{R}^{N_{\vm x}}\times\mathbb{R}^{N_{\vm u}}\times\mathbb{R}^{N_{\vm p}}\times\mathbb{R} \rightarrow \mathbb{R}^{N_{\vm h}}$
as well as for the corresponding terminal constraints
$\vm g_T:\mathbb{R}^{N_{\vm x}}\times\mathbb{R}^{N_{\vm p}}\times\mathbb{R} \rightarrow \mathbb{R}^{N_{\vm g_T}}$
and 
$\vm h_T:\mathbb{R}^{N_{\vm x}}\times\mathbb{R}^{N_{\vm p}}\times\mathbb{R} \rightarrow \mathbb{R}^{N_{\vm h_T}}$, 
respectively.
Finally, \eqref{eq:OCP_boxedUconst} and \eqref{eq:OCP_boxedpTconst}
represent box constraints for the optimization variables $\vm u=\vm
u(t)$, $\vm p$ and $T$ (if applicable).

In comparison to the previous version of the GRAMPC
toolbox~\cite{Kaepernick2014}, the problem
formulation~\eqref{eq:OCP_orig} supports optimization with respect to
parameters, a free end time, general state-dependent equality and inequality
constraints, as well as terminal constraints. Furthermore,
semi-implicit dynamics with a constant mass matrix $\vm M$ can be handled
using the Rosenbrock solver RODAS~\cite{Hairer:Book:1996:Stiff} as numerical integrator. 
This extends the range of possible applications besides MPC to general
optimal control, parameter optimization, and moving horizon
estimation. However, the primary target is embedded model predictive
control of nonlinear systems, as the numerical solution algorithm
is optimized for time and memory efficiency. 
\subsection{Application to model predictive control}
Model predictive control relies on the iterative solution of
an optimal control problem of the form
\begin{subequations}\label{eq:MPC_orig}
\begin{alignat}{4} 
&\min_{\vm u} \quad
&& J(\vm u;\vm x_k) = 
V(\vm x(T)) + \int_0^T l(\vm x(\tau), \vm u(\tau),\tau) \, \dd \tau
\hspace{-4cm} && \hspace{2.5cm}
\label{eq:MPC_cost}
\\
& \textrm{\,s.t.} \quad
&& \vm M \vm{\dot x}(\tau) = \vm f(\vm x(\tau), \vm u(\tau),t_k+\tau) \,,
\quad & \vm x(0) &= \vm x_k
\label{eq:MPC_dynamics}
\\
&&& \vm x(\tau) \in \left[\vm x_{\min}, \vm x_{\max}\right]\,,
\quad & \vm x(T) &\in \Omega_\beta
\label{eq:MPC_const}
\\ \label{eq:MPC_boxedUconst}
&&& \vm u(\tau) \in \left[\vm u_{\min}, \vm u_{\max}\right]
\end{alignat}
\end{subequations}
with the MPC-internal time coordinate $\tau\in[0,T]$ over the prediction
horizon $T$. The initial state value $\vm x_k$ is the measured or
estimated system state at the current sampling instant
$t_k=t_0+k \Delta t$, $k \in \mathbb{N}$ with sampling time $0<\Delta
t\le T$.
The first part of the computed control trajectory $\vm u(\tau)$,
$\tau\in[0,\Delta t)$ is used as control input 
for the actual plant over the time interval $t\in[t_k,t_{k+1})$,
before OCP~\eqref{eq:MPC_orig} is solved again with the new
initial state $\vm x_{k+1}$.

A popular choice of the cost functional~\eqref{eq:MPC_cost}
is the quadratic form
\begin{equation}
V(\vm x) = \|\vm x-\vm x_\text{des}\|^2_{\vm P}\,, \quad
l(\vm x,\vm u)= \|\vm x-\vm x_\text{des}\|_{\vm Q}^2 + 
\|\vm u-\vm u_\text{des}\|_{\vm R}^2
\end{equation}
with the desired setpoint 
$(\vm x_\text{des},\vm u_\text{des})$ and the 
positive (semi-)definite matrices $\vm P$, $\vm Q$, $\vm R$.
Stability is often ensured in MPC by imposing a
terminal constraint $\vm x(T)\in\Omega_\beta$, where the 
set $\Omega_\beta=\{\vm x\in\mathbb{R}^{N_{\vm x}}\,|\, V(\vm x)\le\beta\}$
for some $\beta>0$ is defined in terms of the terminal cost $V(\vm x)$
that can be computed from solving a Lyapunov or Riccati equation that
renders the set $\Omega_\beta$ invariant under a local feedback
law~\cite{Chen:AUT:1998,Mayne:AUT:2000}. 
In view of the OCP formulation \eqref{eq:OCP_orig}, the
terminal region as well as general box constraints on the state as 
given in \eqref{eq:MPC_const} can be expressed as
\begin{equation}
\vm h(\vm x) = \begin{bmatrix}
\vm x-\vm x_\text{max} \\ \vm x_\text{min}-\vm x
\end{bmatrix} \le \vm 0 \,, \quad
\vm h_T(\vm x) = V(\vm x) - \beta \le \vm 0\,.
\end{equation}
Note, however, that terminal constraints are often
omitted in embedded or real-time MPC in order to minimize the
computational effort~\cite{Limon:TAC:2006,Graichen2010,Gruene:AUT:2013}. 
In particular, real-time feasibility is typically achieved by limiting
the number of iterations per sampling step and using the current
solution for warm starting in the next MPC
step, in order to incrementally reduce the suboptimality over the
runtime of the MPC~\cite{Diehl:IEE:2004,Graichen2010,Graichen:AUT:2012}.

An alternative to the ``classical'' MPC
formulation~\eqref{eq:MPC_orig} is shrinking horizon MPC, see
e.g.~\cite{Diehl:SIAM:2005,Skaf:IJRNC:2010,Gruene:MTNS:2014}, where
the horizon length $T$ is shortened over the MPC steps. 
This can be achieved by formulating
the underlying OCP \eqref{eq:MPC_orig} as a free end time problem with a terminal
constraint $\vm g_T(\vm x(T))=\vm x(T)-\vm x_\text{des}=\vm 0$
to ensure that a desired setpoint $\vm x_\text{des}$ is reached
in finite time instead of the asymptotic behavior of
fixed-horizon MPC. 

\subsection{Application to moving horizon estimation} \label{sec:MHE}

Moving horizon estimation (MHE) can be seen as the dual of MPC for
state estimation problems. Similar to MPC, MHE relies on the online solution of
a dynamic optimization problem of the form
\begin{subequations}\label{eq:MHE_1_orig}
\begin{alignat}{3}
& \min_{\vm{\hat x}_k} \quad && 
J(\vm{\hat x}_k; \vm u, \vm y) = \int_{t_k-T}^{t_k} \|
\vm{\hat y}(t) - \vm y(t) \|^2 \, \dd t
\label{eq:MHE_1_cost}
\\
& \textrm{\,s.t.} &&
\vm M\vm{\dot{\hat x}}(t) = \vm f(\vm{\hat x}(t), \vm u(t),t) \,,\quad
\vm{\hat x}(t_k) = \vm{\hat x}_k 
\label{eq:MHE_1_dynamics}
\\ &&& \label{eq:MHE_1_y}
\vm{\hat y}(t) = \vm\sigma(\vm{\hat x}(t))
\end{alignat}
\end{subequations}
that depends on the history of the control $\vm u(t)$ and measured output 
$\vm y(t)$ over the past time window $[t_k-T, t_k]$.
The solution of \eqref{eq:MHE_1_orig} yields the estimate $\vm{\hat x}_k$ 
of the current state $\vm x_k$ such that the estimated output function
\eqref{eq:MHE_1_y} best matches the measured
output $\vm y(t)$ over the past horizon $T$.
Further constraints can be added to the formulation of
\eqref{eq:MHE_1_orig} to incorporate a priori knowledge. 

\GRAMPC can be used for moving horizon estimation by handling the
system state at the beginning of the estimation horizon as
optimization variables, i.e.\ $\vm p=\vm{\hat x} (t_k - T)$.
In addition, a time transformation is required to map $t \in [t_k-T, t_k]$ to
the new time coordinate $\tau \in [0, T]$ along with the corresponding
coordinate transformation  
\begin{align}
\vm{\tilde x}(\tau) = \vm{\hat x}(t_k\!-\!T\!+\!\tau) - \vm p \,,\quad
\vm{\tilde u}(\tau) = \vm u(t_k\!-\!T\!+\!\tau) \,,\quad
\vm{\tilde y}(\tau) = \vm y(t_k\!-\!T\!+\!\tau) \label{eq:coordTrafo}
\end{align}
with the initial condition $\vm{\tilde x}(0) = \vm 0$. 
The optimization problem~\eqref{eq:MHE_1_orig} then
becomes
\begin{subequations}\label{eq:MHE_orig}
\begin{alignat}{2}
\min_{\vm p} \quad
& J(\vm p; \vm{\tilde u}, \vm{\tilde y}) = \int_0^T \| \vm{\hat
  y}(\tau) - \vm{\tilde y}(\tau) \|^2 \, \dd \tau
\\
\textrm{s.t.} \quad
& \vm M\vm{\dot{\tilde x}}(\tau) = \vm f(\vm{\tilde x}(\tau) + \vm p, \vm{\tilde u}(\tau), t_k\!-\!T\!+\!\tau) \,,\quad 
\vm{\tilde x}(0) = \vm 0 
\\ & \vm{\hat y}(\tau) = \vm \sigma(\vm{\tilde x}(\tau)+\vm p)\,.
\end{alignat}
\end{subequations}
The solution $\vm p=\vm{\hat x} (t_k - T)$
of~\eqref{eq:MHE_orig} and the coordinate transformation 
\eqref{eq:coordTrafo} are used to compute the current state estimate with
\begin{equation}
	\vm{\hat x}_k = \vm p + \vm{\tilde x}(T)\,,
\end{equation}
where $\vm{\tilde x}(T)$ is the end point of the state trajectory
returned by \GRAMPC.
In the next sampling step, the parameters can be re-initialized with
the predicted estimate $\vm{\hat x}(t_k-T+\Delta t) = \vm p + \vm{\tilde x}(\Delta t)$.

  Note that the above time transformation can alternatively be reversed
  in order to directly estimate the current state
  $\vm p=\vm{\hat x}(t_k)$. This, however, requires the reverse time
  integration of the dynamics, which is numerically unstable if the
  system is stable in forward time. Vice versa, a reverse
  time transformation is to be preferred for MHE of an unstable process.
\section{Optimization algorithm}
\label{sec:OptAlg}
The optimization algorithm underlying \GRAMPC uses an augmented
Lagrangian formulation in combination with a real-time projected gradient
method. Though SQP or interior point methods are typically superior in
terms of convergence speed and accuracy, the augmented Lagrangian
framework is able to rapidly provide a suboptimal solution at low
computational costs, which is important in view of real-time
applications and embedded optimization.
In the following, the augmented Lagrangian formulation and the
corresponding optimization algorithm are described for solving 
OCP~\eqref{eq:OCP_orig}. 
The algorithm follows a first-optimize-then-discretize approach in
order to maintain the dynamical system structure in the optimality
conditions, before numerical integration is applied.
\subsection{Augmented Lagrangian formulation}
\label{sec:ALFormulation}
The basic idea of augmented Lagrangian methods is to replace the
original optimization problem by its dual problem, see
for example~\cite{Bertsekas1996,Nocedal2006,Boyd2004} as 
well
as~\cite{Fortin:book:1983,Ito:MP:1990,Bergounioux:JOTA:1997,Aguiar:CCA:2016}
for corresponding approaches in optimal control and function space settings.

The augmented Lagrangian formulation adjoins the
constraints~\eqref{eq:OCP_eqconstr}, \eqref{eq:OCP_ieqconstr}
to the cost functional~\eqref{eq:OCP_cost} by means of multipliers 
$\vm{\bar \lag}=(\vm\lag_{\vm g}, \vm\lag_{\vm h}, \vm\lag_{\vm g_T}, \vm\lag_{\vm h_T})$ 
and additional quadratic penalty terms with the penalty parameters 
$\vm{\bar \pen}=(\vm \pen_{\vm g}, \vm \pen_{\vm h}, \vm \pen_{\vm g_T}, \vm \pen_{\vm h_T})$.
A standard approach in augmented Lagrangian theory is to transform the 
inequalities~\eqref{eq:OCP_ieqconstr} into equality constraints by means of
slack variables, which can be analytically solved for~\cite{Bertsekas1996}.
This leads to the overall set of equality constraints (see
Appendix~\ref{sec:AppIeqc} for details)
\begin{subequations}
\label{eq:ieqcII}
\begin{equation}
\label{eq:eqconstr_stacked_2}
\vm{\bar g}(\vm x, \vm u, \vm p, t, \vm{\lag_h},\vm{\pen_h}) = \begin{bmatrix}
\vm g(\vm x, \vm u, \vm p, t) \\
\vm{\bar h}(\vm x, \vm u, \vm p, t, \vm{\lag_h}, \vm{\pen_h})
\end{bmatrix} = \vm 0
\end{equation}
\begin{equation}
\vm{\bar g}_T(\vm x, \vm p, T, \vm \lag_{\vm h_T}, \vm \pen_{\vm h_T}) = \begin{bmatrix}
\vm g_T(\vm x, \vm p, T) \\
\vm{\bar h}_T(\vm x, \vm p, T, \vm \lag_{\vm h_T}, \vm \pen_{\vm h_T})
\end{bmatrix} = \vm 0
\end{equation}
\end{subequations}
with the transformed inequalities
\begin{subequations}\label{eq:transfieqc}
\begin{equation}
\label{eq:transfieqc_1}
\vm {\bar h}(\vm x, \vm u, \vm p, t, \vm \lag_{\vm h}, \vm \pen_{\vm h}) =
\vm\max\left\{ \vm h(\vm x, \vm u, \vm p, t), -\vm \penmatrix_{\vm h}^{-1} \vm \lag_{\vm h} \right\}
\end{equation}
\begin{equation}
\label{eq:transfieqc_2}
\vm {\bar h}_T(\vm x, \vm p, T, \vm \lag_{\vm h_T}, \vm \pen_{\vm h_T}) = 
\vm\max\left\{ \vm h_T(\vm x, \vm p, T), -\vm \penmatrix_{\vm h_T}^{-1} \vm \lag_{\vm h_T} \right\}
\end{equation}
\end{subequations}
and the
diagonal matrix syntax $\vm \penmatrix=\diag(\vm \pen)$. 
The vector-valued $\vm\max$-function is to be understood component-wise.
The equalities \eqref{eq:eqconstr_stacked_2} are adjoined to
the cost functional
\begin{align}
\label{eq:cost_augm_1}
\bar J(\vm u, \vm p, T, \vm{\bar \lag}, \vm{\bar \pen};\vm x_0) = 
\bar V(\vm x, \vm p, T, \vm \lag_T, \vm \pen_T) 
+ \int_0^T \bar l(\vm x, \vm u, \vm p, t, \vm \lag, \vm \pen) \, \mathrm dt
\end{align}
with 
the augmented terminal and integral cost terms
\begin{subequations}
\label{eq:costIII_V&L}
\begin{multline}
\bar V(\vm x, \vm p, T, \vm \lag_T, \vm \pen_T)  = \\
V(\vm x, \vm p, T) 
+ \vm\lag_T\TT \vm{\bar g}_T(\vm x, \vm p, T, \vm \lag_{\vm h_T},\vm \pen_{\vm h_T})
+ \frac12 \| \vm{\bar g}_T(\vm x, \vm p, T, \vm \lag_{\vm h_T},\vm \pen_{\vm h_T}) \|^2_{\vm \penmatrix_T}
\end{multline}
\begin{multline}
\bar l(\vm x, \vm u, \vm p, t, \vm \lag, \vm \pen) = \\
l(\vm x, \vm u, \vm p, t) 
+ \vm\lag\TT \vm{\bar g}(\vm x, \vm u, \vm p, t, \vm{\lag_h},\vm{\pen_h})
+ \frac12 \| \vm{\bar g}(\vm x, \vm u, \vm p, t, \vm{\lag_h},\vm{\pen_h}) \|^2_{\vm C}
\end{multline}
\end{subequations}
and the stacked penalty and multiplier vectors 
$\vm\lag_T = [\vm\lag_{\vm g_T}\TT,\vm\lag_{\vm h_T}\TT]\TT$,
$\vm\pen_T = [\vm\pen_{\vm g_T}\TT,\vm\pen_{\vm h_T}\TT]\TT$,
and
$\vm\lag = [\vm\lag_{\vm g}\TT,\vm\lag_{\vm h}\TT]\TT$,
$\vm\pen = [\vm\pen_{\vm g}\TT,\vm\pen_{\vm h}\TT]\TT$, respectively.
The augmented cost functional~\eqref{eq:cost_augm_1} allows one to
formulate the max-min-problem
\begin{subequations}\label{eq:OCPIII}
\begin{align}
\max_{\vm{\bar\lag}} \, \min_{\vm u, \vm p, T} \quad& 
\costIII J(\vm u, \vm p, T, \vm{\bar\lag}, \vm{\bar\pen};\vm x_0) 
\\ \label{eq:OCPIII:dyn}
\textrm{s.t.} \quad& \vm M \vm{\dot x}(t) = \vm f(\vm x, \vm u, \vm p, t) 
\,,\quad 
\vm x(0) = \vm x_0
\\
& \vm u(t) \in [\vm u_{\min}, \vm u_{\max}]
\,,\quad 
t\in[0,T]
\label{eq:OCPIII:uconst}\\
&\vm p \in [\vm p_{\min}, \vm p_{\max}]
\,,\quad 
T \in [T_{\min}, T_{\max}]\,. \label{eq:OCPIII:ptconst}
\end{align}
\end{subequations}
Note that the multipliers $\vm\lag = [\vm\lag_{\vm g}\TT,\vm\lag_{\vm h}\TT]\TT$
corresponding to the constraints \eqref{eq:OCP_ieqconstr},
respectively~\eqref{eq:transfieqc}, 
are functions of time $t$. In the implementation of
\GRAMPC, the corresponding penalties $\vm\pen_{\vm g}$ and 
$\vm\pen_{\vm h}$ are handled time-dependently as well.

If strong duality holds and $(\vm u^*,\vm p^*,T^*)$ and
$\vm{\bar\lag}^*$ are primal and dual optimal points, they 
form a saddle-point in the sense of
\begin{multline}\label{eq:saddle}
\costIII J(\opt{\vm u},\opt{\vm p},\opt{T}, \vm{\bar\lag}, \vm{\bar\pen};\vm x_0) \leq 
\costIII J(\opt{\vm u},\opt{\vm p},\opt{T}, \opt{\vm{\bar\lag}},\vm{\bar\pen};\vm x_0) \leq 
\costIII J({\vm u},{\vm p},{T}, \opt{\vm{\bar\lag}},\vm{\bar\pen};\vm x_0) \,.
\end{multline} 
On the other hand, if \eqref{eq:saddle} is satisfied, then $(\vm u^*,\vm p^*,T^*)$ and
$\vm{\bar\lag}^*$ are primal and dual optimal and strong
duality holds~\cite{Boyd2004,Allaire:book:2007}.
Moreover, if the saddle-point condition is satisfied for the unaugmented
Lagrangian, i.e.\ for $\vm{\bar\pen}=\vm 0$, then it
holds for all $\vm{\bar\pen}>\vm 0$ and vice versa~\cite{Fortin:book:1983}.

Strong duality typically relies on convexity, which is difficult to
investigate for general constrained nonlinear optimization problems. 
However, the augmented Lagrangian formulation is favorable in this
regard, as the duality gap that may occur for unpenalized, nonconvex
Lagrangian formulations can potentially be closed by the augmented
Lagrangian formulation~\cite{Rockafellar:SIAM:1974}. 

The motivation behind the algorithm presented in
the following lines 
is to solve the dual problem~\eqref{eq:OCPIII} instead of the original
one~\eqref{eq:OCP_orig} by approaching the
saddle-point~\eqref{eq:saddle} from both sides.
In essence, the max-min-problem~\eqref{eq:OCPIII} is solved in an
alternating manner by performing the inner minimization with a projected
gradient method and the outer maximization via a steepest ascent
approach. 
Note that the dynamics~\eqref{eq:OCPIII:dyn} are captured inside the minimization
problem instead of treating the dynamics as equality constraints of
the form $\vm M \vm{\dot x} - \vm f(\vm x, \vm u, \vm p, t) = \vm 0$
in the augmented Lagrangian~\cite{Hager:SIAM:1990}. 
This ensures the dynamical consistency of the computed
trajectories in each iteration of the algorithm, which is important
for an embedded, possibly suboptimal implementation.
\subsection{Structure of the augmented Lagrangian algorithm}
\label{sec:alg_max}
The basic iteration structure of the augmented Lagrangian algorithm is
summarized in Algorithm~1 and will be detailed in the
Sections~\ref{sec:alg_min} to \ref{sec:update_mult_pen}.
The initialization of the algorithm concerns the multipliers
$\vm{\bar\lag}^1$ and penalties $\vm{\bar\pen}^1$ as well as the definition of
several tolerance values that are used for the convergence check (Section~\ref{sec:convergence})
and the update of the multipliers and penalties in \eqref{eq:alg_mult} and
\eqref{eq:alg_pen}, respectively.

In the current augmented Lagrangian iteration $i$, 
the inner minimization is carried out by solving the
OCP~\eqref{eq:OCP_minProb} for the current
set of multipliers $\vm{\bar\lag}^i$ and penalties $\vm{\bar\pen}^i$.
Since the only remaining constraints within~\eqref{eq:OCP_minProb} are
box constraints on the optimization variables~\eqref{eq:OCPIV:uconst}  
and \eqref{eq:OCPIV:ptconst}, the problem can be efficiently solved by the
projected gradient method described in Section~\ref{sec:alg_min}.
The solution of the minimization step consists of the 
control vector $\vm u\ind {i}$ and of the parameters $\vm p\ind {i}$ and free end time 
$T\ind {i}$, if these are specified in the problem at hand.
The subsequent convergence check in Algorithm~1
rates the constraint violation as well as convergence behavior of the
previous minimization step and is detailed in Section~\ref{sec:convergence}.

If convergence is not reached yet, the multipliers and penalties are
updated for the next iteration of the algorithm, as detailed in
Section~\ref{sec:update_mult_pen}. 
Note that the penalty update in~\eqref{eq:alg_pen} relies on the last
two iterates of the constraint functions~\eqref{eq:alg_constr}.
In the initial iteration $i=1$ and if \GRAMPC is used within an MPC setting, the constraint 
functions $\vm g^0$, $\vm h^{0}$, $\vm g_T^{0}$, $\vm h_T^{0}$ are
warm-started by the corresponding last iterations of the previous MPC
run. Otherwise, the penalty update is started in iteration $i=2$.
\vspace{3mm}
{\parindent0mm
\myalg{1}{Augmented Lagrangian algorithm}

\footnotesize
        \vspace{1mm}
        \textbf{Initialization}
        
        \begin{itemize}\itemsep1mm
        \item Initialize multipliers $\vm{\bar\lag}^0$ and penalties $\vm{\bar\pen}^0$
        \item Set tolerances 
          ($\vm{\varepsilon}_{\vm{g}}>0$, $\vm{\varepsilon}_{\vm{h}}>0$, $\vm{\varepsilon}_{\vm{g}_T}>0$, 
          $\vm{\varepsilon}_{\vm{h}_T}>0$, $\varepsilon_\text{rel,c}>0$)
        \end{itemize}
        
        \vspace{-3mm}\hrulefill\vspace{0.5mm}

        \parindent0mm
        \textbf{for} $i=1$ \textbf{to} $i_\text{max}$ \textbf{do}
        
        \vspace{-1mm}
	\begin{itemize}\itemsep1mm		
        \item Compute $(\vm u^i,\vm p^i,T^i)$ 
          with $\vm u^i\!=\!\vm u^i(t)$ 
          and state $\vm x^i\!=\!\vm x^i(t)$, $t\!\in\![0,T^i]$ by solving
          the optimal control problem \text{(see Section~\ref{sec:alg_min})}
          \begin{subequations}\label{eq:OCP_minProb}
            \begin{align}
              \min_{\vm u, \vm p, T} \quad
              & \costIII J(\vm u, \vm p, T, \vm{\bar\lag}\ind{i}, \vm{\bar\pen}\ind{i};\vm x_0) 
                \label{eq:OCPIV_costfct} \\
              \textrm{\!s.t.} \quad
              & \vm M \vm{\dot x}(t) = \vm f(\vm x, \vm u, \vm p, t) \,,\quad 	\vm x(0) = \vm x_0	
                \label{eq:OCPIV_dynamics} \\
              & \vm u(t) \in [\vm u_{\min}, \vm u_{\max}]\,,
                \quad t\in[0,T]                          
                \label{eq:OCPIV:uconst} \\
              & \vm p \in [\vm p_{\min}, \vm p_{\max}]
                \,,\quad T \in [T_{\min}, T_{\max}]
                \label{eq:OCPIV:ptconst}
            \end{align}
          \end{subequations}
          
        \item Store constraint functions 
          
          \vspace{-4mm}
          \begin{subequations}
            \label{eq:alg_constr}
            \begin{alignat}{3}
              \vm g^i(t) &= \vm g(\vm x^i(t), \vm u^i(t),t)\,,
              \quad & 
              \vm g_T^i &= \vm g_T(\vm x^i(T),T^i)
              \\
              \vm{\bar h}^i(t) &=\vm{\bar h}(\vm x^i(t), \vm
              u^i(t),\vm p^i,t,\vm\lag_{\vm h}^i(t),\vm\pen_{\vm h}^i(t)) \,,
              \quad &
              \vm h_T^i &= \vm{\bar h}_T(\vm x^i(T),\vm p^i,T^i,\vm\lag_{\vm h_T}^i,\vm\pen_{\vm h_T}^i)
            \end{alignat}
          \end{subequations}
          \vspace{-4mm}

        \item If convergence criterion is
          reached (see Section~\ref{sec:convergence}),
          \textbf{break}
		
        \item Update multipliers $\vm{\bar\lag}^i=(\vm\lag_{\vm g}^i,
          \vm\lag_{\vm h}^i,\vm\lag_{\vm g_T}^i,\vm\lag_{\vm h_T}^i)$
          according to (see Section~\ref{sec:update_mult_pen}) 
          
          \vspace{-4mm}
          \begin{subequations}\label{eq:alg_mult}
          \begin{alignat}{3} \label{eq:alg_mult_g}
              \hspace{-3mm}
              \vm{\lag}\ind{i+1}_{\vm{g}}(t) &=
              \vm{\zeta}_{\vm{g}}\!\left(\vm{\lag}\ind{i}_{\vm{g}}(t),\vm{\pen}\ind{i}_{\vm{g}}(t),
                \vm{g}\ind{i}(t),\vm{\varepsilon}_{\vm{g}}\right),
              \quad & 
              \vm{\lag}\ind{i+1}_{\vm{g}_T} &=
              \vm{\zeta}_{\vm{g}_T}\!\left(\vm{\lag}\ind{i}_{\vm{g_T}},\vm{\pen}\ind{i}_{\vm{g_T}},
                \vm{g}_T\ind{i},\vm{\varepsilon}_{\vm{g}_T}\right)
              \\[-1mm] \label{eq:alg_mult_h}
              \vm{\lag}\ind{i+1}_{\vm{h}}(t) &= \vm{\zeta}_{\vm{h}}\!\left(\vm{\lag}\ind{i}_{\vm{h}}(t), 
                \vm{\pen}\ind{i}_{\vm{h}}(t),\vm{\bar h}\ind{i}(t),\vm{\varepsilon}_{\vm{h}}\right), 
              \quad &
              \vm{\lag}\ind{i+1}_{\vm{h}_T} &= \vm{\zeta}_{\vm{h}_T}\!\left(
                \vm{\lag}\ind{i}_{\vm{h}_T},\vm{\pen}\ind{i}_{\vm{h}_T},\vm{\bar h}_T\ind{i},
                \vm{\varepsilon}_{\vm{h}_T}\right)
            \end{alignat}
          \end{subequations}
          	
        \vspace{-3mm}
        \item Update penalties $\vm{\bar\pen}^i=(\vm\pen_{\vm g}^i,
          \vm\pen_{\vm h}^i,\vm\pen_{\vm g_T}^i,\vm\pen_{\vm h_T}^i)$, if
          $i\ge 2$ or MPC warmstart 
                (see Section~\ref{sec:update_mult_pen}) 

                \vspace{-4mm}
                \begin{subequations}\label{eq:alg_pen}
		\begin{alignat}{3}\label{eq:alg_pen_g}
                  \vm{\pen}\ind{i+1}_{\vm{g}}(t) &=
                  \vm{\xi}_{\vm{g}}\!\left(\vm{\pen}\ind{i}_{\vm{g}}(t),\vm{g}\ind{i}(t),
                    \vm{g}\ind{i-1}(t),\vm{\varepsilon}_{\vm{g}}\right),
                  \quad &
                  \vm{\pen}\ind{i+1}_{\vm{g}_T} &= \vm{\xi}_{\vm{g}_T}\!\left(
                    \vm{\pen}\ind{i}_{\vm{g_T}},\vm{g}_T\ind{i},\vm{g}_T\ind{i-1},
                    \vm{\varepsilon}_{\vm{g}_T}\right)                  
                  \\[-1mm] \label{eq:alg_pen_h}
                  \vm{\pen}\ind{i+1}_{\vm{h}}(t) &= 
                  \vm{\xi}_{\vm{h}}\!\left(\vm{\pen}\ind{i}_{\vm{h}}(t),\bar{\vm h}\ind{i}(t),
                    \bar{\vm h}\ind{i-1}(t),\vm{\varepsilon}_{\vm{h}}\right),
                  \quad &
                  \vm{\pen}\ind{i+1}_{\vm{h}_T} &= \vm{\xi}_{\vm{h}_T}\!\left(
                    \vm{\pen}\ind{i}_{\vm{h}_T},\bar{\vm h}_T\ind{i},
                    \bar{\vm h}_T\ind{i-1},\vm{\varepsilon}_{\vm{h}_T}\right)
                      \end{alignat}
                    \end{subequations}
                \end{itemize}

                \vspace{-3mm}
                \textbf{end}
\vspace{2mm}
\hrule
}

\vspace{5mm}
Note that if the end time $T$ is treated as optimization variable as
shown in Algorithm~1, the evaluation of~\eqref{eq:alg_pen}
would formally require to redefine the constraint functions
$\vm g^{i-1}(t)$ and $\vm h^{i-1}(t)$, $t\in[0,T^{i-1}]$ from the previous
iterations to the new horizon length $T^{i}$ by either shrinkage or extension. 
This redefinition is not explicitly stated in Algorithm~1, since
the actual implementation of \GRAMPC stores the trajectories in
discretized form, which implies that only the discretized time
vector must be recomputed, once the end time $T^i$ is updated.

\subsection{Gradient algorithm for inner minimization problem}
\label{sec:alg_min}

The OCP~\eqref{eq:OCP_minProb} inside the augmented Lagrangian
algorithm corresponds to the inner minimization problem of the
max-min-formulation~\eqref{eq:OCPIII} for the current iterates of the
multipliers $\vm{\bar\lag}^i$ and $\vm{\bar\pen}^i$. 
A projected gradient method is used to solve
OCP~\eqref{eq:OCP_minProb} to a desired accuracy or for a fixed number
of iterations.

The gradient algorithm relies on the solution of the first-order
optimality conditions defined in terms of the Hamiltonian
\begin{equation}\label{eq:Hamiltonian}
H(\vm x, \vm u,\vm p,\vm \lambda, t, \vm\mu,\vm c) = 
\costIII{l}(\vm x,\vm u, \vm p,t,\vm\mu,\vm c)+\vm\lambda\TT\vm f(\vm x, \vm u,\vm p, t)
\end{equation}
with the adjoint states $\vm \lambda \in \mathbb{R}^{N_{\vm x}}$.
In particular, the gradient algorithm iteratively solves the canonical
equations, see e.g.~\cite{Cao:SIAM:2003},
\begin{subequations}\label{eq:optcond}
  \begin{alignat}{3}
    \vm M\vm{\dot x} &= \vm f (\vm x,\vm u,\vm p,t)\,,
    \quad &\opt{\vm x}(0) & = \vm{x}_0, \label{eq:optxcond}
    \\
    \vm M\TT \vm{\dot \lambda} &= -
    H_{\vm{x}}(\vm x,\vm u,\vm p,\vm\lambda,t, \vm\mu,\vm c)\,,
    \quad  &\vm M\TT \vm\lambda(T) &=	
    \costIII{V}_{\vm x}(\vm x(T),\vm p,T,\vm\mu_T,\vm c_T)
    \label{eq:optlcond}
  \end{alignat}
\end{subequations}
consisting of the original dynamics~\eqref{eq:optxcond} and the
adjoint dynamics~\eqref{eq:optlcond} in forward and backward time
and computes a gradient update for the control in order to minimize
the Hamiltonian in correspondence with Pontryagin's Maximum Principle 
\cite{Kirk1970,Berkovitz1974}, i.e.\
\begin{equation}
\label{eq:optinputcond}
\underset{{\vm{u} \in [\vm u_{\min}, \vm u_{\max}]}}{\min} 
H\big(\vm x(t),\vm u,\vm p,\vm\lambda(t),t,\vm\mu(t),\vm c(t)\big) \,, \quad
t\in[0,T]\,.
\end{equation}
If parameters $\vm p$ and/or the end time $T$ are additional
optimization variables, the corresponding gradients have to be
computed as well. 
Algorithm~2 lists the overall projected gradient
algorithm for the full optimization case, i.e.\ for the optimization
variables $(\vm u,\vm p,T)$, for the sake of completeness.

\newpage
\vspace{3mm}
{\parindent0mm
\myalg{2}{Projected gradient algorithm (inner minimization)}

\footnotesize
  \vspace{1mm}
  \textbf{Initialization}
        
  \vspace{-1mm}
  \begin{itemize}\itemsep0.5mm
  \item Initialize $\vm u^{i|1}(t)$, $t\in[0,T^{i|1}]$ and $\vm p^{i|1}$, $T^{i|1}$
  \item Compute $\vm x^{i|1}(t)$, $t\in[0,T^{i|1}]$ by solving~\eqref{eq:alg:fwdint}
    for $j=0$
  \item Set step size adaptation factors $\gamma_{\vm p}>0$, $\gamma_{T}>0$
  \end{itemize}

  \vspace{-3mm}\hrulefill\vspace{0.5mm}
        
  \parindent0mm
  \textbf{for} $j=1$ \textbf{to} $j_\text{max}$ \textbf{do}
  
  \begin{itemize}\itemsep1mm		

  \item Compute $\vm\lambda^{i|j}(t)$ by backward integration of
    \begin{multline}
      \label{eq:alg2:adj}
      \vm M\TT \vm{\dot\lambda}^{i|j}(t) = 
      - H_{\vm x}\!\left(\vm x^{i|j}(t),\vm u^{i|j}(t),\vm p^{i|j}, 
      \vm\lambda^{i|j}(t),t, \vm\mu^i(t),\vm c^i(t) \right), \\
      \vm M\TT \vm\lambda^{i|j}(T^{i|j}) = \bar V_{\vm x}\!\left(\vm
        x^{i|j}(T^{i|j}), T^{i|j},\vm\mu^{i}_T, \vm c_T^{i}\right)
    \end{multline}
    
  \item Compute gradients 
    with 
    $H^{i|j}(t):=H\big({\vm x}\ind {i|j}(t),{\vm u}\ind{i|j}(t),
        {\vm p}\ind {i|j},{\vm \lambda}\ind {i|j}(t),t, \vm\mu^i(t),\vm c^i(t)\big)$ 
    \begin{subequations}\label{eq:ls_dir}
      \begin{align}
        \vm d_{\vm u} \ind {i|j}(t) &= 
        H_{\vm u}\big({\vm x}\ind {i|j}(t),{\vm u}\ind{i|j}(t),
        {\vm p}\ind {i|j},{\vm \lambda}\ind {i|j}(t),t, \vm\mu^i(t),\vm c^i(t)\big)
        \label{eq:ls_diru} \\[-1mm]
        \vm d_{\vm p} \ind {i|j}~ &=
      \bar V_{\vm p}(\vm x^{i|j}(T^{i|j}),\vm p^{i|j}, T^{i|j},\vm\mu_T^{i},\vm c_T^{i})+
      \int_0^{T^{i|j}} \!\! H_{\vm p}^{i|j}(t) \,\dd t  \label{eq:ls_dirp}\\
      d_T \ind {i|j}~ &= 
      \bar V_T(\vm x^{i|j}(T^{i|j}),\vm p^{i|j}, T^{i|j},\vm\mu_T^{i},\vm c_T^{i})+
      H^{i|j}(T^{i|j})
      \label{eq:ls_dirT}
      \end{align}
    \end{subequations}

  \item Compute step size $\alpha\ind {i|j}$ by solving the line search problem
\begin{equation}\label{eq:ls_prob}
\underset{\alpha>0}{\min} ~~
\costIII{J}\Big(
\vm\psi_{\vm u}\!\left(\vm u\ind{i|j} - \alpha\,\vm  d_{\vm u} \ind {i|j}\right)\!,
\vm \psi_{\vm p}\!\left(\vm p \ind {i|j} - \gamma_{\vm p}\,\alpha\,\vm  d_{\vm p} \ind {i|j}\right)\!,
\psi_T\!\left(T\ind{i|j} - \gamma_{T}\,\alpha\,\vm  d_T \ind {i|j} \right);\vm x_0 \Big)
\end{equation}

\item Update $(\vm u^{i|j+1},\vm p^{i|j+1},T^{i|j+1})$ according to
\begin{subequations}
\label{eq:alg2:update}
\begin{alignat}{2}
\vm u \ind {i|j+1}(t) &= \vm\psi_{\vm u}\!\left(\vm u \ind {i|j}(t) - 
\alpha\ind {i|j} \vm  d_{\vm u} \ind {i|j}(t)\right)
\\
\vm p \ind {i|j+1} &= \vm \psi_{\vm p}\!\left(\vm p \ind {i|j} - 
\gamma_{\vm p}\,\alpha\ind {i|j} \vm  d_{\vm p} \ind {i|j}\right)
\\
T \ind {i|j+1} &= \psi_T\!\left(T\ind{i|j} - \gamma_{T}\,\alpha\ind {i|j} d_T \ind {i|j}\right)
\end{alignat}
\end{subequations}

\item Compute $\vm x^{i|j+1}(t)$ by forward integration of
  \begin{equation}
    \label{eq:alg:fwdint}
    \vm M \vm{\dot x}^{i|j+1}(t) = \vm f\!\left(\vm x^{i|j+1}(t),\vm u^{i|j+1}(t),
      \vm p^{i|j+1},t\right), \quad \vm x^{i|j+1}(0) = \vm x_0
  \end{equation}
  
\item Evaluate convergence criterion
\begin{equation}\label{eq:minimizationconv}
\eta\ind{i|j+1} = 
\max\left\{\frac{\|\vm u \ind {i|j+1}-\vm u \ind {i|j}\|_{L_2}} {\|\vm u \ind {i|j+1}\|_{L_2}} \,,
			\frac{\zweinorm{\vm p \ind {i|j+1}-\vm p \ind {i|j}}}{\zweinorm{\vm p \ind {i|j+1}}}\,,
			 \frac{|T\ind {i|j+1}-T\ind {i|j}|}{T\ind {i|j+1}} \right\} 
\end{equation}

\item If $\eta\ind{i|j+1} \leq \varepsilon_\text{rel,c}$ or
  $j=j_\text{max}$, \textbf{break}

\end{itemize}

\vspace{-1mm}
\textbf{end}

\vspace{-1mm}\hrulefill\vspace{0.5mm}

\textbf{Output to Algorithm~1}

\vspace{-1mm}
\begin{itemize}\itemsep1mm
  \item Set 
    $\vm u^i(t) := \vm u^{i|j+1}(t)$ and
    $\vm x^i(t) := \vm x^{i|j+1}(t)$ 
  \item Set $\vm p \ind {i} := \vm p \ind {i|j+1}$ and
    $T^i := T^{i|j+1}$
    \item Return $\eta^i:=\eta^{i|j+1}$
\end{itemize}

\vspace{-1mm}
\hrule
}\vspace{3mm}

\newpage
The gradient algorithm is initialized with an initial control 
$\vm u^{i|1}(t)$ and initial parameters $\vm p^{i|1}$ and time length
$T^{i|1}$. In case of MPC, these initial values are taken from the
last sampling step using a warmstart strategy with an optional time
shift in order to compensate for the horizon shift by the sampling
time $\Delta t$.

The algorithm starts in iteration $j=1$ with computing the
corresponding state trajectory $\vm x^{i|j}(t)$ as well as the adjoint
state trajectory $\vm\lambda^{i|j}(t)$ by integrating the adjoint
dynamics~\eqref{eq:alg2:adj} in reverse time.
In the next step, the gradients \eqref{eq:ls_dir} are computed 
and the step size $\alpha^{i|j}>0$ is determined from the
line search problem \eqref{eq:ls_prob}. 
The projection functions $\vm \psi_{\vm u}(\vm u)$, $\vm \psi_{\vm p}(\vm p)$, 
and $\psi_{T}(T)$ project the inputs, the parameters, and the end time 
onto the feasible sets \eqref{eq:OCPIV:uconst} and
\eqref{eq:OCPIV:ptconst}. For instance, 
$\vm \psi_{\vm u}(\vm u)=[\psi_{\vm u,1}(u_1) \ldots \psi_{\vm u,N_{\vm u}}(u_{N_{\vm u}})]\TT$ is
defined by
\begin{equation}
\psi_{\vm u,i}(u_i) = 
\begin{cases} u_i & \text{if~} u_i \in (u_{\text{min},i},u_{\text{max},i})
\\
u_{\text{min},i} & \text{if~} u_i \le u_{\text{min},i} \,,
\\
u_{\text{max},i} & \text{if~} u_i \ge u_{\text{max},i}
\end{cases}
\quad
i=1,\ldots,N_{\vm u} \,.
\end{equation}
The next steps in Algorithm~2 are the
updates~\eqref{eq:alg2:update} of the control $\vm u^{i|j+1}$,
parameters $\vm p^{i|j+1}$, and end time $T^{i|j+1}$ as well as the
update of the state trajectory $\vm x^{i|j+1}$ in~\eqref{eq:alg:fwdint}. 

The convergence measure $\eta^{i|j+1}$ in \eqref{eq:minimizationconv}  
rates the relative gradient changes of $\vm u$, $\vm p$, and $T$.
If the gradient scheme converges in the sense of 
$\eta\ind{i|j+1} \leq \varepsilon_\text{rel,c}$ with  
threshold $\varepsilon_\text{rel,c}>0$ or if the maximum number of
iterations $j_\text{max}$ is reached,  
the algorithm terminates and returns the last solution to
Algorithm~1.
Otherwise, $j$ is incremented and the gradient iteration continues.

An important component of Algorithm~2 is the line
search problem \eqref{eq:ls_prob}, which is performed in all search
directions simultaneously. The scaling factors $\gamma_{\vm p}$ and
$\gamma_{T}$ can be used to scale the step sizes relative to each
other, if they are not of the same order of magnitude or if the
parameter or end time optimization is highly sensitive.  
\GRAMPC implements two different line search strategies, an adaptive and an
explicit one, in order to solve \eqref{eq:ls_prob} in an accurate and
robust manner without involving too much computational load.

The adaptive strategy evaluates the cost
functional~\eqref{eq:ls_prob} for three different step sizes, i.e.\
$(\alpha_i,\bar J_i)$, $i=1,2,3$ with $\alpha_1<\alpha_2<\alpha_3$,
in order to compute a polynomial fitting function $\Phi(\alpha)$ of the form
\begin{equation}\label{eq:ls_adapt_approx}
\Phi(\alpha) = p_0+p_1\alpha +p_2\alpha^2\,,
\end{equation}
where the constants $p_i$ are computed from the test
points $(\alpha_i,\bar J_i)$, $i=1,2,3$. The approximate step size can
then be analytically derived by solving
\begin{equation}
\alpha \ind {i|j} = \underset{\alpha \in [\alpha_1,\alpha_3]}{\arg\min} \,\,\,\Phi(\alpha)\,.
\end{equation}
The interval $[\alpha_1,\alpha_3]$ is adapted in the next gradient
iteration, if the step size $\alpha\ind j$ is close to the interval's borders, 
see \cite{Kaepernick2014} for more details.

Depending on the OCP or MPC problem at hand, the adaptive line search
method may not be suited for time-critical applications, since the
approximation of the cost function~\eqref{eq:ls_adapt_approx} requires
to integrate the system dynamics~\eqref{eq:OCPIV_dynamics}
and the cost integral~\eqref{eq:OCPIV_costfct} three times. 
This computational load can be further reduced by an explicit line
search strategy, originally discussed in \cite{Barzilai1988} and
adapted to the optimal control case in \cite{Kaepernick2013}.
Motivated by the secant equation in quasi-Newton
methods~\cite{Barzilai1988}, this strategy minimizes the
difference between two updates of the optimization variables 
$(\vm u,\vm p,T)$ for the same step size and without considering the
corresponding box constraints \eqref{eq:OCPIV:uconst} and
\eqref{eq:OCPIV:ptconst}. The explicit method solves the problem 
\begin{alignat}{2} \label{eq:ls_expl_prob}
& \underset{\alpha > 0}{\min}~ \bigl\| {\vm u}\ind{i|j+1} - {\vm u}\ind {i|j} \bigr\|^2_{L_2}
+ \bigl\|{\vm p}\ind{i|j+1} - {\vm p}\ind {i|j} \bigr\|^2_2
+ \big({T}\ind{i|j+1} - {T}\ind {i|j} \big)^2 
= \\ \nonumber
& \underset{\alpha>0}{\min}~
\bigl\| \Delta\vm u^{i|j} -\alpha\Delta \vm d^{i|j}_{\vm u} \bigr\|^2_{L_2} 
+
\bigl\| \Delta\vm p^{i|j} -\gamma_{\vm p}\alpha \Delta \vm d^{i|j}_{\vm p} \bigr\|^2_2
+
\bigl\| \Delta T^{i|j} -\gamma_{T}\alpha \Delta d^{i|j}_T \bigr\|^2_2 \,.
\end{alignat}
where $\Delta$ denotes the difference between the last and current
iterate, e.g.\ $\Delta\vm u^{i|j}=\vm u^{i|j}-\vm u^{i|j-1}$. The
analytic solution is given by
\begin{equation}
\label{eq:ls_expl1}
\alpha\ind{i|j} = 
\frac{\langle \Delta \vm{u}\ind{i|j},\Delta\vm{d}\ind{i|j}_{\vm u}\rangle 
	+\gamma_{\vm p} \langle \Delta\vm{p}\ind{i|j},\Delta\vm{d}\ind{i|j}_{\vm p}\rangle 
	+\gamma_{T} \Delta T\ind{i|j}\Delta{d}\ind{i|j}_{T}}
{\langle \Delta\vm{d}\ind{i|j}_{\vm u},\Delta \vm{d}\ind{i|j}_{\vm u} \rangle 
	+ \gamma_{\vm p}^2 \langle{\Delta\vm{d}\ind{i|j}_{\vm p}},\Delta \vm{d}\ind{i|j}_{\vm p}\rangle 
	+ \gamma_{T}^2\big( \Delta {d}\ind{i|j}_{T}\big)^2}\,.
\end{equation}
Alternatively, the minimization 
\begin{equation*} 
\underset{\alpha > 0}{\min}~ \bigl\| \alpha\Delta  {\vm u}\ind{i|j} -\Delta \vm{d}\ind{i|j}_{\vm u} \bigr\|^2_{L_2}
+ \bigl\| \gamma_{\vm p}\alpha\Delta {\vm p}\ind{i|j} -\Delta \vm{d}\ind{i|j}_{\vm p} \bigr\|_2^2
+\big(\gamma_{T}\alpha\Delta {T}\ind{i|j} -\Delta {d}\ind{i|j}_{T}\big)^2
\end{equation*}
can be carried out, similar to \eqref{eq:ls_expl_prob}, leading to
the corresponding solution  
\begin{equation}
\label{eq:ls_expl2}
\alpha\ind{i|j} = 
\frac{\langle \Delta \vm{u}\ind{i|j}, \Delta\vm{u}\ind{i|j}\rangle 
	+ \gamma_{\vm p} \langle\Delta\vm{p}\ind{i|j}, \Delta\vm{p}\ind{i|j} \rangle
	+ \gamma_{T} \big(\Delta T\ind{i|j}\big)^2}
{\langle \Delta\vm{u}\ind{i|j}, \Delta \vm{d}\ind{i|j}_{\vm u} \rangle 
	+ \gamma_{\vm p}^2 \langle\Delta\vm{p}\ind{i|j}, \Delta \vm{d}\ind{i|j}_{\vm p} \rangle 
	+ \gamma_{T}^2 \Delta {T}\ind{i|j} \Delta {d}\ind{i|j}_{T}} \,.
\end{equation}
Both explicit formulas \eqref{eq:ls_expl1} and \eqref{eq:ls_expl2} are implemented in
\GRAMPC as an alternative to the adaptive line search strategy.

\subsection{Convergence criterion}
\label{sec:convergence}

The gradient scheme in Algorithm~2 solves the inner
minimization problem~\eqref{eq:OCP_minProb} of the augmented
Lagrangian algorithm and returns the solution 
$(\vm u^i,\vm p^i,T^i)$ as well as the maximum relative gradient
$\eta^i$ that is computed in~\eqref{eq:minimizationconv} and used to
check convergence inside the gradient algorithm.
The outer augmented Lagrangian iteration 
in Algorithm~1 also uses this criterion 
along with the convergence check of the constraints, i.e.\
\begin{equation}
\begin{bmatrix}
|{\vm {g}}\ind{i}_T| \\ \vm\max\{ \bar{\vm h}\ind{i}_T, \vm 0 \}
\end{bmatrix}
\le
\begin{bmatrix}
\vm{\varepsilon_{\vm{g}_T}} \\ \vm{\varepsilon_{\vm{h}_T}} 
\end{bmatrix}
\,\,\,\land\,\,\,
\underset{t\in[0,T]}{\vm\max}
\begin{bmatrix}
|{\vm {g}}\ind{i}(t)|
\\[1mm]
\vm\max\{ \vm{\bar h}^i(t), \vm 0 \}
\end{bmatrix}
\le
\begin{bmatrix}
\vm{\varepsilon_{\vm{g}}} \\ \vm{\varepsilon_{\vm{h}}} 
\end{bmatrix}
\,\,\,\land\,\,\, 
\eta\ind{i} \leq \varepsilon_\text{rel,c} \,.
\end{equation}
The thresholds
$\vm{\varepsilon_{\vm{g}}},\vm{\varepsilon_{\vm{g_T}}},
\vm{\varepsilon_{\vm{h}}}, \vm{\varepsilon_{\vm{h_T}}}$ are
vector-valued to rate each constraint individually.
If the maximum number of augmented Lagrangian iterations is reached, 
i.e. $i=i_\text{max}$, the algorithm terminates in order to ensure
real-time feasibility. Otherwise, $i$ is incremented and the next
minimization of problem~\eqref{eq:OCP_minProb} is carried out.

\subsection{Update of multipliers and penalties}
\label{sec:update_mult_pen}

The multiplier update~\eqref{eq:alg_mult} in Algorithm~1
is carried out via a steepest ascent approach, whereas the penalty
update~\eqref{eq:alg_pen} uses an adaptation strategy that rates
the progress of the last two iterates. 
In more detail, the multiplier update function for a 
single equality constraint $g^i:=g(\vm x^i,\vm u^i,\vm p^i,t)$
is defined by
	\begin{equation}\label{eq:lag_update_eqc}
		\zeta_{g}( \lag_g^i,\pen_g^i,g^i,\varepsilon_g) 
		= 
		\begin{cases}
			\lag_g^i + (1-\rho)\pen_g^i g^i
			& \text{if } \left|g^i\right| >  \varepsilon_g
			\,\land \,\eta^i \leq \varepsilon_\text{rel,u}
			\\
			\lag_g^i
			& \text{else}\,.
		\end{cases}
	\end{equation}
The steepest ascent direction is the residual of the constraint $g^i$
with the penalty $\pen_g^i$ serving as step size parameter.
The additional damping factor $0\leq\rho\leq1$ is introduced to
increase the robustness of the multiplier update in the augmented
Lagrangian scheme.
In the \GRAMPC implementation, the 
multipliers are additionally limited by an upper bound $\lag_\text{max}$,  
in order to avoid unlimited growth and numerical stability issues.
The update \eqref{eq:lag_update_eqc} is skipped if the constraint is
satisfied within the tolerance $\varepsilon_g$ or if the gradient method is
not sufficiently converged, which is checked by the maximum relative
gradient $\eta^i$, cf.\ \eqref{eq:minimizationconv}, and the update
threshold $\varepsilon_\text{rel,u}>0$. \GRAMPC uses a larger value
for $\varepsilon_\text{rel,u}$ than the 
convergence threshold $\varepsilon_\text{rel,c}$ of the
gradient scheme in Algorithm~2. This accounts for the
case that the gradient algorithm might not have converged to the desired
tolerance before the maximum number of iterations $j_\text{max}$ is
reached, e.g.\ in real-time MPC applications, where only one or two 
gradient iterations might be applied.
In this case, $\varepsilon_\text{rel,u}\gg\varepsilon_\text{rel,c}$
ensures that the multiplier update is still performed provided that
at least ``some'' convergence was reached by the inner minimization.

The penalty $\pen_g^i$ corresponding to the equality constraint $g^i$ is
updated according to the heuristic update function
\begin{equation}
  \label{eq:pen_update_eqc}
  \xi_g(\pen_g^i,g^i,g^{i-1}\!,\varepsilon_g)=
  \begin{cases}
    \beta_{\mathrm{in}}\pen_g^i & \text{if }
    \big|g^i\big| \!\geq \max\big\{\gamma_{\mathrm{in}}\big|g^{i-1}\big|
    , \varepsilon_g\big\}
        \land \eta\ind{i} \!\leq \varepsilon_\text{rel,u}
    \\ 
    \beta_{\mathrm{de}} \,\pen_g^i & \text{else if } 
    \left|g^i\right| \leq \gamma_{\mathrm{de}}\,\varepsilon_g
    \\ 
    \pen_g^i & \text{else}
  \end{cases}
\end{equation}
that is basically motivated by the LANCELOT package \cite{Conn2013,Nocedal2006}.  
The penalty $\pen_g^i$ is increased by the factor $\beta_\text{in}\ge 1$,
if the last gradient scheme converged, i.e. $\eta\ind{i}\le\varepsilon_\text{rel,u}$, 
but insufficient progress (rated by $\gamma_\text{in}>0$) was made by the constraint
violation compared to the previous iteration $i-1$. 
The penalty is decreased by the factor $\beta_\text{de}\le 1$ if the
constraint $g^i$ is sufficiently satisfied within its tolerance with
$0<\gamma_\text{de}<1$. The setting $\beta_\text{de}=\beta_\text{in}=1$
can be used to keep $\pen_g^i$ constant.
Similar to the multiplier update~\eqref{eq:lag_update_eqc}, 
\GRAMPC restricts the penalty to upper and lower bounds
$\pen_\text{min}\le \pen_g^i \le \pen_\text{max}$, in order to avoid
negligible values as well as unlimited growth of $\pen_g^i$. 
The lower penalty bound $\pen_\text{min}$ is particularly relevant in case
of MPC applications, where typically only a few iterations are
performed in each MPC step. \GRAMPC provides a support function that
computes an estimate of $\pen_\text{min}$ for the MPC problem at hand.

The updates \eqref{eq:lag_update_eqc} and \eqref{eq:pen_update_eqc} 
define the vector functions $\vm\zeta_{\vm g}$, $\vm\zeta_{\vm g_T}$
and $\vm\xi_{\vm g}$, $\vm\xi_{\vm g_T}$ in \eqref{eq:alg_mult_g} and
\eqref{eq:alg_pen_g} with $N_{\vm g}$ and $N_{\vm g_T}$ components,
corresponding to the number of equality and terminal equality
constraints. 
Note that the multipliers for the equality and inequality constraints are
time-dependent, i.e.\ $\vm\lag_{\vm g}(t)$ and $\vm\lag_{\vm h}(t)$, which
implies that the functions $\vm\zeta_{\vm g}$ and $\vm\xi_{\vm g}$,
resp.~\eqref{eq:lag_update_eqc} and \eqref{eq:pen_update_eqc}, are
evaluated pointwise in time.

The inequality constrained case is handled in a
similar spirit. For a single inequality constraint 
$\bar h^i:=\bar h(\vm x^i,\vm u^i,\vm p^i,t,\lag_h,\pen_h)$, cf.\ 
\eqref{eq:transfieqc_1}, the multiplier and penalty updates are
defined by
\begin{equation}
  \label{eq:lag_update_ieqc}
  \zeta_{h}( \lag_h^i,\pen_h^i, \bar h^i,\varepsilon_h) 
  = 
  \begin{cases}
    \lag_h^i +  (1-\rho) \pen_h^i \bar h^i
    & \text{if } \big(\bar h^i > \varepsilon_h
    \land \eta^i \leq \varepsilon_\text{rel,u}
    \big) \lor \bar h^i < 0
    \\
    \lag_h^i
    & \text{else}
  \end{cases}
\end{equation}
and 
\begin{equation}
  \label{eq:pen_update_ieqc}
  \xi_h(c_h^i,\bar h^i,\bar h^{i-1},\varepsilon_h) =
  \begin{cases}
    \beta_{\mathrm{in}} \,c_h^i
    & \text{if }
    \bar h^i \geq \max\big\{\gamma_{\mathrm{in}} \bar h^{i-1}, \varepsilon_h\big\}
    \land \eta\ind{i} \leq \varepsilon_\text{rel,u}
    \\ 
    \beta_{\mathrm{de}} \, c_h^i
    & \text{else if } \bar h^i \leq \gamma_{\mathrm{de}} \, \varepsilon_h
    \\ 
    c_h^i  & \text{else}
    \end{cases}
\end{equation}
and constitute the vector functions 
$\vm\zeta_{\vm h}$, $\vm\zeta_{\vm h_T}$
and $\vm\xi_{\vm h}$, $\vm\xi_{\vm h_T}$ in \eqref{eq:alg_mult_h} and
\eqref{eq:alg_pen_h} with $N_{\vm h}$ and $N_{\vm h_T}$ components. 
The condition $\bar h^i < 0$ in \eqref{eq:lag_update_ieqc} ensures
that the Lagrangian multiplier $\lag_h^i$ is reduced for inactive
constraints, which corresponds to either $h^i < 0$ or
$-\frac{\lag_h}{\pen_h}$ in view of the transformation
\eqref{eq:transfieqc}. 
\section{Structure and usage of GRAMPC}
\label{sec:StrUsage}

This section describes the framework of \GRAMPC and 
illustrates its general usage. \GRAMPC is designed to be portable and
executable on different operating systems and hardware without the use of
external libraries. 
The code is implemented in plain C with a user-friendly interface to C++, \Matlab/\Simulink, and \Dspace. 
The following lines give an overview of \GRAMPC and demonstrate how to
implement and solve a problem. 

\subsection{General structure}
\label{sec:StrUsage:Str}

\begin{figure}[t]
  \begin{center}
  	  \includegraphics[scale=0.85]{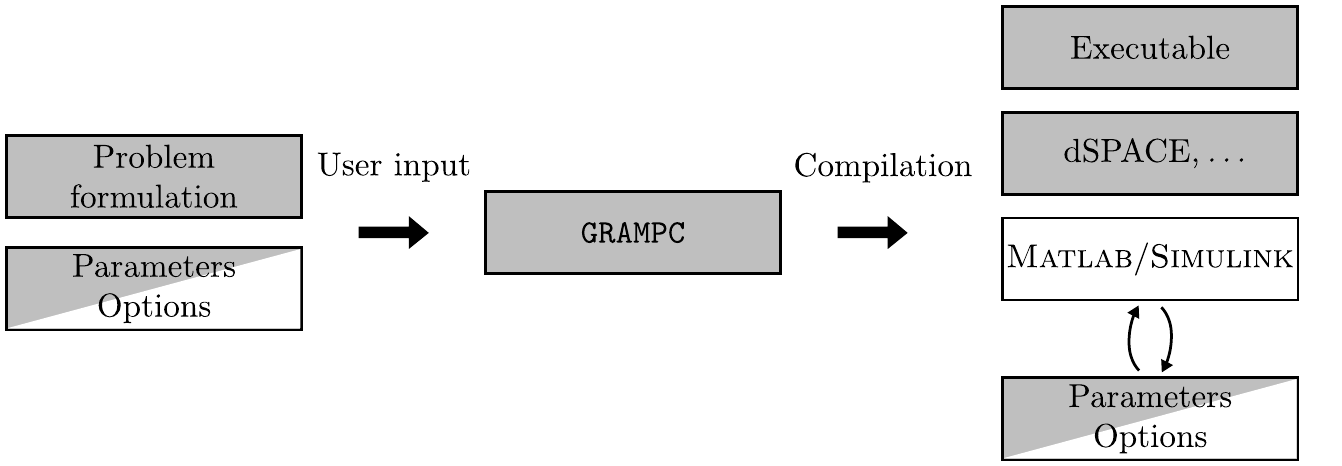}
      \caption{General structure of \GRAMPC (grey: C code, white:
        \Matlab level).}
      \label{fig:StrUsage:StructureGRAMPC}
      \vspace{-3mm}
  \end{center}
\end{figure}

\Fig{\ref{fig:StrUsage:StructureGRAMPC}} shows the general structure
of \GRAMPC and the steps that are necessary to compile an executable
\GRAMPC project. 
The first step in creating a new project is to define the problem
using the provided C function templates, which will be detailed more in
\Sec{\ref{sec:StrUsage:PrblDef}}. 
The user has the possibility to set problem specific parameters and
algorithmic options concerning the numerical integrations in the
gradient algorithm, the line search strategy as well as further
preferences, also see \Sec{\ref{sec:StrUsage:Opt}}.

\begin{figure}[t]
  \begin{center}
  	  \includegraphics[scale=1.0]{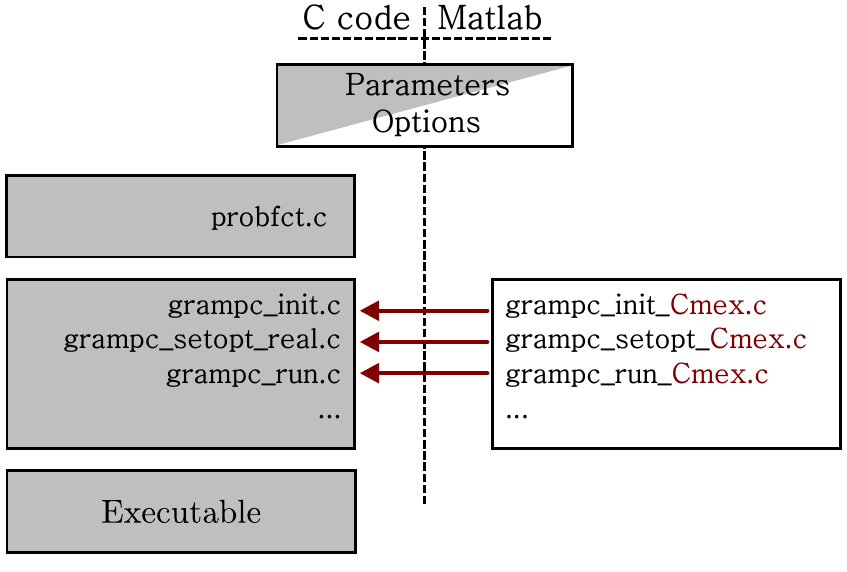}
      \caption{Interfacing of \GRAMPC to C (grey) and Matlab
        (white). Each \GRAMPC function written in plain C has a
        corresponding Cmex function.}
      \label{fig:StrUsage:MexInterface}
      \vspace{-5mm}
  \end{center}
\end{figure}

A specific problem can be parameterized and numerically
solved using C/C++, \Matlab/\Simulink, or \Dspace. 
A closer look on this functionality is given in
Figure~\ref{fig:StrUsage:MexInterface}. 
The workspace of a \GRAMPC project as well as algorithmic options and
parameters are stored by the structure variable \verb+grampc+. Several
parameter settings are problem-specific and need to be provided,
whereas other values are set to default values. 
A generic interface allows one to manipulate the \verb+grampc+
structure in order to set algorithmic options or parameters for the problem at hand. 
The functionalities of \GRAMPC can be manipulated from
\Matlab/\Simulink by means of mex routines that are wrappers for the
corresponding C functions. 
This allow one to run a \GRAMPC project
with different parameters and options without recompilation. 

\subsection{Problem definition}
\label{sec:StrUsage:PrblDef}
The problem formulation in \GRAMPC follows a generic structure.
The essential steps for a problem definition are illustrated 
for the following MPC example 
 \begin{subequations}
 \label{eq:struct:ballonplate}
\begin{alignat}{2}
\label{eq:struct:ballonplate:cost}
& \min_{u(\cdot)} & \quad &J(u;\vm x_k) = 
\frac{1}{2}\Delta\vm x\TT(T) \vm P \Delta \vm x(T)+
\frac{1}{2}\int_{0}^{T} \Delta \vm x\TT(T)
\vm Q
\Delta\vm x +
R 
\Delta u^2 \,\dd \tau \\
\label{eq:struct:ballonplate:dynamics}
& \text{\;s.t.} & \quad & 
\begin{bmatrix} \dot x_1 \\ \dot x_2\end{bmatrix}
=
\begin{bmatrix} 0 & 1 \\ 0 & 0\end{bmatrix}\begin{bmatrix} x_1\\ x_2\end{bmatrix}
+
\begin{bmatrix} -0.04 \\ -7.01\end{bmatrix}u
\,, \quad
\begin{bmatrix} x_1(0)\\ x_2(0)\end{bmatrix} = 
\begin{bmatrix} x_{k,1}\\ x_{k,2}\end{bmatrix}
\\
\label{eq:struct:ballonplate:stateconstr}
&& \quad & 
\begin{bmatrix} -0.2\\-0.1\end{bmatrix}
\le\begin{bmatrix} x_1\\ x_2\end{bmatrix}\le
\begin{bmatrix} 0.01\\0.1\end{bmatrix}
\,, \quad
|u| \le 0.0524
\end{alignat}
with $\Delta \vm x=\vm x-\vm x_\text{des}$, $\Delta u=u-u_\text{des}$,
and the weights
\begin{equation}
\vm P=\vm Q = \begin{bmatrix} 100 & 0 \\ 0 & 10\end{bmatrix}, \quad
R = 1
\,.
\end{equation}
\end{subequations}
The dynamics~\eqref{eq:struct:ballonplate:dynamics} are a simplified
linear model of a single axis of a ball-on-plate system~\cite{Richter2012} that is
also included in the testbench of \GRAMPC (see Section~\ref{sec:performance}).
The horizon length and the sampling time are set to $T=0.3$\,s and $\Delta t=10\,$ms, respectively.

The problem formulation \eqref{eq:struct:ballonplate} is provided in
the user template \verb+probfct.c+. The following listing gives an
expression of the function structure and the problem implementation
for the ball-on-plate example~\eqref{eq:struct:ballonplate}. 
 
\begin{lstlisting}
/* Problem dimensions */
void ocp_dim(typeInt *Nx, typeInt *Nu, typeInt *Np, typeInt *Ng, typeInt 
             *Nh, typeInt *NgT, typeInt *NhT, typeUSERPARAM *userparam)
{ 
  *Nx  = 2; *Nu  = 1;  *Np  = 0;
  *Nh  = 4; *Ng  = 0;
  *NgT = 0; *NhT = 0;
}

/* Right-hand side of dynamics M dx/dt = f(t,x,u,p,userparam) */
void ffct(typeRNum *out, ctypeRNum t, ctypeRNum *x, ctypeRNum *u, 
          ctypeRNum *p, typeUSERPARAM *userparam)
{ 
  out[0] = x[1]-0.04*u[0];
  out[1] =     -7.01*u[0];
}

/* Integral cost function int l(t,x(t),u(t),p,...) dt */
void lfct(typeRNum *out, ctypeRNum t, ctypeRNum *x, ctypeRNum *u, ctypeRNum
          *p, ctypeRNum *xdes, ctypeRNum *udes, typeUSERPARAM *userparam)
{
  typeRNum* param = (typeRNum*)userparam;
  out[0] = 0.5*(param[0]*(x[0]-xdes[0])*(x[0]-xdes[0]) + \
                param[1]*(x[1]-xdes[1])*(x[1]-xdes[1]) + \
                param[2]*(u[0]-udes[0])*(u[0]-udes[0]) );
}

/* Terminal cost function V(T,x(T),p,xdes,userparam) */
void Vfct(typeRNum *out, ctypeRNum T, ctypeRNum *x, ctypeRNum *p, 
          ctypeRNum *xdes, typeUSERPARAM *userparam)
{
  typeRNum* param = (typeRNum*)userparam;
  out[0] = 0.5*(param[3]*(x[0]-xdes[0])*(x[0]-xdes[0]) + \
                param[4]*(x[1]-xdes[1])*(x[1]-xdes[1]) );
}

/* Inequality constraints h(t,x(t),u(t),p,uperparam) <= 0 */
void hfct(typeRNum *out, ctypeRNum t, ctypeRNum *x, ctypeRNum *u, 
          ctypeRNum *p, typeUSERPARAM *userparam)
{
  typeRNum* param = (typeRNum*)userparam;
  out[0] =  param[5] - x[0];
  out[1] = -param[6] + x[0];
  out[2] =  param[7] - x[1];
  out[3] = -param[8] + x[1];
}

\end{lstlisting}

The naming of the functions follows the nomenclature of the general
OCP formulation~\eqref{eq:OCP_orig}, except for the function
\verb+ocp_dim+, which defines the dimensions of the optimization
problem. Problem specific parameters can 
be used inside the single functions via the pointer \verb+userparam+.
For convenience, the desired setpoint (\verb+xdes+,\verb+udes+) to be
stabilized by the MPC is provided to the cost functions separately and
therefore does not need to be passed via \verb+userparam+.

In addition to the single OCP functions, \GRAMPC requires the
derivatives w.r.t.\ state $\vm x$, control $\vm u$, and if applicable
w.r.t.\ the optimization parameters $\vm p$ and end time $T$, in order to evaluate the
optimality conditions in Algorithm~2. 
Jacobians that occur in multiplied form, see e.g.\ 
$(\frac{\partial \vm f}{\partial \vm x}\big)\TT\vm \lambda$ 
in the adjoint dynamics~\eqref{eq:alg2:adj}, have to be provided in this
representation. This helps to avoid unnecessary zero multiplications 
in case of sparse Jacobians. The following listing shows an excerpt of
the corresponding gradient functions.

\begin{lstlisting}
/* Multiplied Jacobian (df/dx)^T * mult */
void dfdx_mult(typeRNum *out, ctypeRNum t, ctypeRNum *x, ctypeRNum *mult, 
               ctypeRNum *u, ctypeRNum *p, typeUSERPARAM *userparam)
{
  out[0] = 0;
  out[1] = mult[0];
}
...
/* Jacobian dl/dx */
void dldx(typeRNum *out, ctypeRNum t, ctypeRNum *x, ctypeRNum *u, ctypeRNum 
          *p, ctypeRNum *xdes, ctypeRNum *udes, typeUSERPARAM *userparam)
{
  typeRNum* param = (typeRNum*)userparam;
  out[0] = param[0]*(x[0]-xdes[0]);
  out[1] = param[1]*(x[1]-xdes[1]);
}
...
/* Multiplied Jacobian (dh/dx)^T * mult */
void dhdx_mult(typeRNum *out, ctypeRNum t, ctypeRNum *x, ctypeRNum *u, 
               ctypeRNum *p, ctypeRNum *mult, typeUSERPARAM *userparam)
{
  out[0] = -mult[0]+mult[1];
  out[1] = -mult[2]+mult[3];
}
...
\end{lstlisting}

\subsection{Calling procedure and options}
\label{sec:StrUsage:Opt}

\GRAMPC provides several key functions that are required for
initializing and calling the MPC solver. As shown in
Figure~\ref{fig:StrUsage:MexInterface}, there exist mex wrapper
functions that ensure that the interface for interacting with \GRAMPC
is largely the same under C/C++ and \Matlab.

The following listing gives an idea on how to initialize \GRAMPC and
how to run a simple MPC loop for the ball-on-plate example under \Matlab.

\begin{lstlisting}
  % user parameters
  userparam = [100,10,1,100,10,-0.2,0.01,-0.1,0.1]

  % initialization
  grampc = grampc_init_Cmex(userparam);

  % set parameters
  grampc = grampc_setparam_Cmex(grampc,'x0',[0.1;0.01]);
  grampc = grampc_setparam_Cmex(grampc,'xdes',[-0.2;0]);
  grampc = grampc_setparam_Cmex(grampc,'u0',0);
  grampc = grampc_setparam_Cmex(grampc,'udes',0);
  grampc = grampc_setparam_Cmex(grampc,'Thor',0.3);
  grampc = grampc_setparam_Cmex(grampc,'dt',0.01);
  grampc = grampc_setparam_Cmex(grampc,'t0',0);
  grampc = grampc_setparam_Cmex(grampc,'umin',0.0524);
  grampc = grampc_setparam_Cmex(grampc,'umax',-0.0524);

  % set options
  grampc = grampc_setopt_Cmex(grampc,'Nhor',20);
  grampc = grampc_setopt_Cmex(grampc,'MaxGradIter',2);
  grampc = grampc_setopt_Cmex(grampc,'MaxMultIter',3);
  grampc = grampc_setopt_Cmex(grampc,'InequalityConstraints','on');
  grampc = grampc_setopt_Cmex(grampc,'Integrator','heun');

  % MPC loop
  for i = 1:iMPC
    % run GRAMPC
    grampc = grampc_run_Cmex(grampc);
    ...
    % set new initial state
    grampc = grampc_setparam_Cmex(grampc,'x0',grampc.sol.xnext);
    ...
  end 
  ...
\end{lstlisting}

The listing 
also shows some of the algorithmic settings, 
e.g.\ the number of discretization points \verb+Nhor+ for the horizon $[0,T]$, the
maximum number of iterations $(i_\text{max},j_\text{max})$ 
for Algorithm~1 and 2, 
or the choice of integration scheme for solving the canonical
equations~\eqref{eq:alg2:adj}, \eqref{eq:alg:fwdint}.
Currently implemented integration methods are (modified) Euler, Heun,
4th order Runge-Kutta as well as the solver
RODAS \cite{Hairer:Book:1996:Stiff} that implements a 4th order Rosenbrock method for
solving semi-implicit differential-algebraic 
equations with possibly singular and sparse mass matrix $\vm M$, cf.\
the problem definition in \eqref{eq:OCP_orig}.
The Euler and Heun methods use a fixed step size depending on the
number of discretization points (\verb+Nhor+),
whereas RODAS and Runge-Kutta use adaptive step size control.
The choice of integrator therefore has significant impact on the computation
time and allows one to optimize the algorithm in terms of accuracy
and computational efficiency. 
Further options not shown in the listing
are e.g.\ the settings 
(\verb+xScale+, \verb+xOffset+) and (\verb+uScale+,\verb+uOffset+)
in order to scale the input and state variables of the
optimization problem.

The initialization and calling procedure for \GRAMPC is largely the
same under C/C++ and \Matlab. One exception is the handling of user
parameters in \verb+userparam+. Under C, \verb+userparam+ can be an
arbitrary structure, whereas the \Matlab interface restricts
\verb+userparam+ to be of type array (of arbitrary length). 

Moreover, the \Matlab call of \verb+grampc_run_Cmex+ returns an
updated copy of the \verb+grampc+ structure as output argument in
order to comply with the \Matlab policy to not manipulate input arguments.

\section{Performance evaluation}
\label{sec:performance}

The intention of this section is to evaluate the performance of \GRAMPC under realistic conditions and for meaningful problem settings. To this end, an MPC testbench suite is defined to evaluate the computational performance in comparison to other state-of-the-art MPC solvers and to demonstrate the portability of \GRAMPC to real-time and embedded hardware. 
The remainder of this section demonstrates the  
handling of typical problems from the field of MPC, moving horizon estimation and optimal control.

\subsection{General MPC evaluation}

The MPC performance of \GRAMPC is evaluated for a testbench that covers a wide range of meaningful MPC applications. For the sake of comparison, the two MPC toolboxes ACADO and VIATOC are included in the evaluation, although it is not the intention of this section to rigorously rate the performance against other solvers, as such a comparison is difficult to carry out objectively. The evaluation should rather give a general impression about the performance and properties of \GRAMPC.
In addition, implementation details are presented for running the MPC testbench examples with \GRAMPC on \Dspace and ECU level.

\subsubsection{Benchmarks}\label{sec:Benchmarks}

Table~\ref{tab:BenchmarkOverview} gives an overview of the considered MPC benchmark problems in terms of the system dimension, the type of constraints (control/state/general nonlinear constraints), 
the dynamics (linear/nonlinear and explicit/semi-implicit) as well as the respective references.
The MPC examples are evaluated with \GRAMPC as well as with ACADO Toolkit~\cite{Houska2011} and VIATOC~\cite{Kalmari2015}. 

The testbench includes three linear problems (mass-spring-damper, helicopter, ball-on-plate) and one semi-implicit reactor example, where the mass matrix $\vm M$ in the semi-implicit form~\eqref{eq:OCP_dynamics} is computed from the original partial differential equation (PDE) using finite elements, also see Section~\ref{sec:reactor_PDE}.
The nonlinear chain problem is a scalable MPC benchmark with $m$ masses.
Three further MPC examples are defined with nonlinear constraints. The permanent magnet synchronous machine (PMSM) possesses spherical voltage and current constraints in dq-coordinates, whereas the crane example with three degrees of freedom (DOF) and the vehicle problem include a nonlinear  constraint to simulate a geometric restriction that must be bypassed (also see Section~\ref{sec:NC_MPC} for the crane example). Three of the problems (PMSM, 2D-crane, vehicle) are not evaluated with \VIATOC, as nonlinear constraints cannot be handled by \VIATOC at this stage.

\begin{table}[b]
  \caption{Overview of MPC benchmark problems.}
  \label{tab:BenchmarkOverview}
  \vspace{-2mm}
  \centering
  \setlength{\tabcolsep}{3pt}
  \renewcommand*{\arraystretch}{1.1}
  \begin{adjustbox}{max width=\textwidth}
    \begin{tabular}{lccccccccccccccc}
      \toprule
      \multirow{2}{*}{Problem} && \multicolumn{2}{c}{Dimensions} && \multicolumn{3}{c}{Constraints} 
      && \multicolumn{2}{c}{Dynamics} && \multirow{2}{*}{Reference}
      \\
      \cmidrule{3-4} \cmidrule{6-8} \cmidrule{10-11}
      && $N_{\vm x}$ & $N_{\vm u}$ && $\vm u$ & $\vm x$ & nonl. && semi-impl. & linear &&
\\ \midrule
Mass-spring-damper 	&& 10 &  2 && yes &  no &  no &&  no & yes && \cite{Kaepernick:diss:2016} \\
Motor (PMSM)		&&  4 &  2 &&  yes & yes & yes &&  no &  no && \cite{Englert:CEP:2018} \\
Nonl.\ chain ($m\!=\!4$) 	&& 21 &  3 && yes &  no &  no &&  no &  no && \cite{Kirches2012} \\
Nonl.\ chain ($m\!=\!6$) 	&& 33 &  3 && yes &  no &  no &&  no &  no && \cite{Kirches2012} \\
Nonl.\ chain ($m\!=\!8$)	&& 45 &  3 && yes &  no &  no &&  no &  no && \cite{Kirches2012} \\
Nonl.\ chain ($m\!=\!10$) 	&& 57 &  3 && yes &  no &  no &&  no &  no && \cite{Kirches2012} \\
2D-Crane   		&&  6 &  2 && yes & yes & yes &&  no &  no && \cite{Kaepernick2013} \\
3D-Crane 	 	&& 10 &  3 && yes &  no &  no &&  no &  no && \cite{Graichen2010} \\	
Helicopter		&&  6 &  2 && yes & yes &  no &&  no & yes && \cite{Tondel2002} \\	
Quadrotor		&&  9 &  4 && yes &  no &  no &&  no &  no && \cite{Kaepernick2014} \\	
VTOL			&&  6 &  2 && yes &  no &  no &&  no &  no && \cite{Sastry2013} \\
Ball-on-plate		&&  2 &  1 && yes & yes &  no &&  no & yes && \cite{Richter2012} \\
Vehicle			&&  5 &  2 && yes &  no & yes &&  no &  no && \cite{Werling2012} \\
CSTR reactor		&&  4 &  2 && yes &  no &  no &&  no &  no && \cite{Rothfuss1996} \\
PDE reactor		&& 11 &  1 && yes &  no &  no && yes &  no && \cite{Utz2010}
\\
\bottomrule 
\end{tabular}
\end{adjustbox}
\end{table}

For the \GRAMPC implementation, 
most options are set to their default values.
The only adapted parameters concern the horizon length $T$, the number of supporting points for the integration scheme and the integrator itself as well as the
number of augmented Lagrangian and gradient iterations, $i_\text{max}$ and $j_\text{max}$, respectively.
Default settings are used for the multiplier and penalty updates for the sake of consistency, see Algorithm~1 as well as Section~\ref{sec:update_mult_pen}. 
Note, however, that the performance and computation time of \GRAMPC can be further optimized by tuning the parameters related to the penalty update to a specific problem. 
All benchmark problems listed in Table~\ref{tab:BenchmarkOverview} are available as example implementations in \GRAMPC.

\ACADO and \VIATOC are individually tuned for each MPC problem by varying the number of shooting intervals and iterations in order to either achieve 
minimum computation time (setting ``speed'') or optimal performance in terms of minimal cost at reasonable computational load (setting ``optimal''). 
The solution of the quadratic programs of \ACADO was done with qpOASES~\cite{Ferreau:MPC:2014}. 

The single MPC projects are integrated in a closed-loop simulation environment with a fourth-order Runge-Kutta integrator with adaptive step size to ensure an accurate system integration regardless of the integration schemes used internally by the MPC toolboxes.
The evaluation was carried out on a Windows 10 machine with Intel(R) Core(TM) i5-5300U CPU running at $\SI{2.3}{\giga\hertz}$ using the Microsoft Visual C++ 2013 Professional (C) compiler. Each simulation was run multiple times to obtain a reliable average computation time.

\subsubsection{Evaluation results}
\label{sec:runtime_evaluation}

Table~\ref{tab:BenchmarkEvalMPC} shows the evaluation results for the benchmark problems in terms of computation time and cost value integrated over the whole time interval of the simulation scenario. 
The cost values are normalized to the best one of each benchmark problem. 
The results for \ACADO and \VIATOC are listed for the settings ``speed'' and ``optimal'', as mentioned above. 
The depicted computation times are the mean computation times, averaged over the complete time horizon of the simulation scenario. 
The best values for computation time and performance (minimal cost) for each benchmark problem are highlighted in bold font. 

The linear MPC problems (mass-spring-damper, helicopter, ball-on-plate) with quadratic cost functions can be tackled well by \VIATOC and \ACADO due to their underlying linearization techniques. 
The PDE reactor problem contains a stiff system dynamics in semi-implicit form. 
\ACADO can handle such problems well using its efficient integration schemes, whereas \VIATOC relies on fixed step size integrators and therefore requires a relatively large amount of discretization points. 
While \GRAMPC achieves the fastest computation time, the cost value of both \ACADO settings as well as the \VIATOC optimal setting is lower. A similar cost function, however, can be achieved by \GRAMPC when deviating from the default parameters. 

In case of the state constrained 2D-crane problem, the overall cost is higher for \ACADO than for \GRAMPC. This appears to be due to the fact that almost over the complete simulation time a nonlinear constraint of a geometric object to be bypassed is active and \ACADO does not reach the new setpoint in the given time. 
A closer look at this problem is taken in Section~\ref{sec:NC_MPC}.

The CSTR reactor example possesses state and control variables in different orders of magnitude and therefore benefits from scaling. Since \GRAMPC supports scaling natively, the computation time is faster than for \VIATOC, where the scaling would have to be done manually. 
Due to the Hessian approximation used by \ACADO, it is far less affected by the different scales in the states and controls.

A large difference in the cost values occurs for the VTOL example (Vertical Take-Off and Landing Aircraft). Due to the nonlinear dynamics and the corresponding coupling of the control variables, it seems that the gradient method underlying the minimization steps of \GRAMPC is more accurate and robust when starting in an equilibrium position than the iterative linearization steps undertaken by \ACADO and \VIATOC.

\begin{table}[t]
	\caption{Evaluation results for the benchmark problems in Table~\ref{tab:BenchmarkOverview} 
          with overall integrated cost $J_\text{int}$ and computation time $t_\text{CPU}$ in milliseconds.}
	\label{tab:BenchmarkEvalMPC}
        \vspace{-2mm}
	\centering
	\setlength{\tabcolsep}{3pt}
	\renewcommand*{\arraystretch}{1.1}
	\begin{adjustbox}{max width=\textwidth}
		\begin{tabular}{@{}lrrp{0cm}rrp{0cm}rrp{0cm}rrp{0cm}rr@{}}
			\toprule
			\multirow{2}{*}{Problem}  & \multicolumn{2}{c}{\multirow{2}{*}{\GRAMPC}}          && \multicolumn{2}{c}{\ACADO} && \multicolumn{2}{c}{\ACADO} && \multicolumn{2}{c}{\VIATOC} && \multicolumn{2}{c}{\VIATOC} \\
			  & \multicolumn{2}{c}{}          && \multicolumn{2}{c}{(optimal)} && \multicolumn{2}{c}{(speed)} && \multicolumn{2}{c}{(optimal)} && \multicolumn{2}{c}{(speed)} \\
			\cmidrule{2-3} \cmidrule{5-6} \cmidrule{8-9} \cmidrule{11-12} \cmidrule{14-15}
			                                     & $J_\text{int}$~            & $t_\text{CPU}$	   &&  $J_\text{int}$~           & $t_\text{CPU}$     &&  $J_\text{int}$~ & $t_\text{CPU}$     && $J_\text{int}$~     & $t_\text{CPU}$     && $J_\text{int}$~ & $t_\text{CPU}$     \\
			\midrule
			Mass-spring-damper               			  	   & $\mb{1.000}$		    & $0.060$ && $1.030$         & $0.370$ && $1.351$ & $0.110$ && $\mb{1.000}$ & $0.075$ && $1.221$ & $\mb{0.049}$ \\
			Motor (PMSM)  	        & $1.006$ & $\mb{0.032}$ && {$\mb{1.000}$} & $0.129$ && $1.096$ & $0.057$ && ---     & ---     && ---     &  ---	   \\
			Nonl. chain ($m\!=\!4$)   & $1.003$ & $2.022$ && $1.003$ & $4.666$ && $1.027$ & $\mb{1.770}$ && $\mb{1.000}$ & $12.32$ && $1.041$  & $7.660$ \\
			Nonl. chain ($m\!=\!6$)   & $\mb{1.000}$ & $\mb{2.492}$ && $\mb{1.000}$ & $12.438$ && $1.028$ & $5.407$ && $\mb{1.000}$ & $18.10$ && $1.042$ & $11.48$ \\
			Nonl. chain ($m\!=\!8$)   & $\mb{1.000}$ & $\mb{5.155}$ && $1.005$ & $26.834$ && $1.050$ & $10.127$ && $\mb{1.000}$ & $38.51$ && $1.055$ & $15.30$ \\
			Nonl. chain ($m\!=\!10$)  & $\mb{1.000}$ & $\mb{7.810}$ && $1.004$ & $43.467$ && $1.042$ & $21.745$ && $1.004$ & $52.84$ && $1.149$ & $24.85$ \\
			2D-Crane  	& $\mb{1.000}$		& $\mb{0.019}$ && $1.932$        & $0.446$ && $2.110$ & $0.058$ && ---     & ---     && --- & ---      \\
			3D-Crane   & ${1.002}$ 		& $\mb{0.036}$ && $\mbox{1.000}$        & $0.728$ && $1.013$ & $0.321$ && $1.163$ & $0.839$ && $1.160$ & $0.194$ \\
			Helicopter      & $1.017$        & $0.054$ && $1.006$        & $0.163$ && $1.096$ & ${\mb{0.045}}$ && $\mb{1.000}$ & $0.185$ && $1.060$ & $0.071$ \\
			Quadrotor       & $\mb{1.000}$        & $\mb{0.022}$ && $\mb{1.000}$        & $1.465$ && $1.005$ & $0.243$ && $1.009$ & $0.535$ && $1.010$ & $0.113$ \\
			VTOL            & $\mb{1.000}$        & $\mb{0.033}$ && $1.210$         & $0.229$ && $1.227$  & $0.090$ && $1.243$  & $0.092$ && $1.259$  & $0.083$ \\
			Ball-on-plate   & $1.113$         & $\mb{0.014}$ && $1.112$         & $0.068$ && $1.126$  & $0.022$ && $\mb{1.000}$ & $0.116$ && $1.158$ & $0.018$ \\
			Vehicle 		& $1.027$      & $\mb{0.092}$ && $\mb{1.000}$ & $1.313$ && $1.003$ & $0.322$ && --- & --- && --- & --- \\
			CSTR reactor  & $\mb{1.000}$ & $\mb{0.058}$ && $\mb{1.000}$ & $0.195$ && $1.002$ & $0.076$ && $1.016$ & $3.711$ && $1.532$ & $0.176$ \\
			PDE reactor    & $1.059$         & $\mb{0.362}$ && $\mb{1.000}$         & $6.259$ && $1.005$  & $0.498$ && $1.046$  & $7.214$ && $1.072$  & $2.868$ \\
			\bottomrule 
		\end{tabular}
	\end{adjustbox}
\end{table}

\begin{figure}[t]
	\centering
	\includegraphics[page=1]{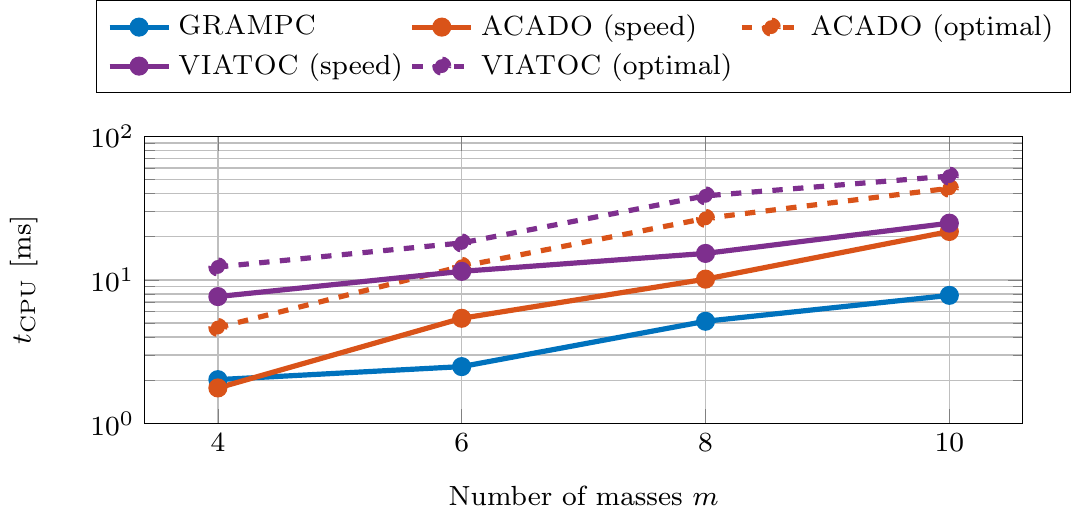}
	\caption{Computation time for the nonlinear chain example (also see Table~\ref{tab:BenchmarkEvalMPC}).}
	\label{fig:NLChainCompTime}
\end{figure}

The scaling behavior of the MPC schemes w.r.t.\ the problem dimension is investigated for the nonlinear chain in Table~\ref{tab:BenchmarkEvalMPC}. Four different numbers of masses are considered, corresponding to 21-57 state variables and three controls. Although the algorithmic complexity of the augmented Lagrangian/gradient projection algorithm of \GRAMPC grows linearly with the state dimension, this is not exactly the case for the nonlinear chain, as the stiffness of the dynamics increases for a larger number of masses, which leads to smaller step sizes of the adaptive Runge-Kutta integrator that was used in \GRAMPC for this problem.
\ACADO shows a more significant increase in computation time for larger values of $m$, which was to be expected in view of the SQP algorithm underlying \ACADO.
The computation time for \VIATOC is worse for this example, since only fixed step size integrators are available in the current release, which requires to increase the number of discretization points manually. \Fig{\ref{fig:NLChainCompTime}} shows a logarithmic plot of the computation time for all three MPC solvers plotted over the number of masses of the nonlinear chain.

\begin{figure}[t]
	\centering
	\includegraphics[page=1]{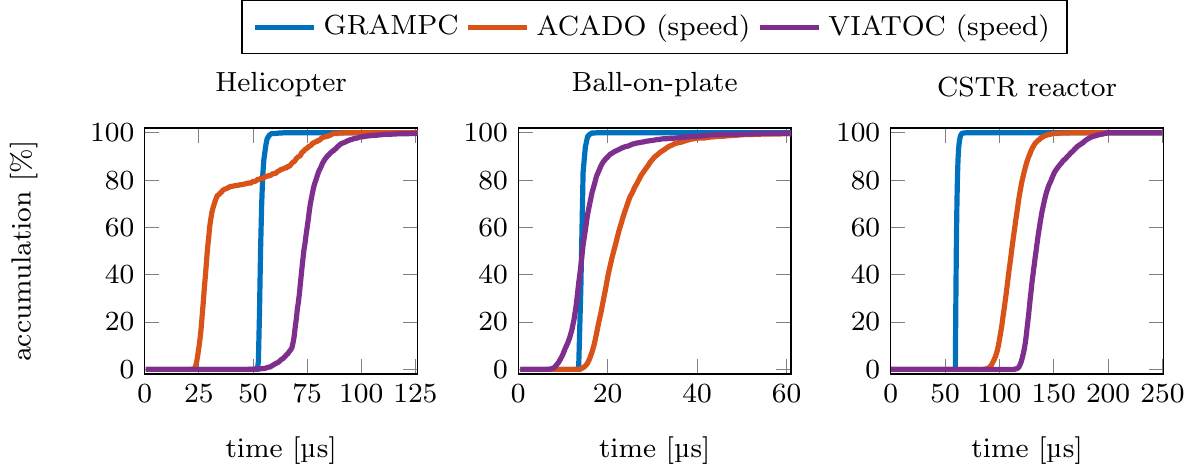}
	\caption{Accumulation of the computation times for three different examples (averaged over 30 runs).}
	\label{fig:compTimeCumSum}
\end{figure}

The computation times shown in Table~\ref{tab:BenchmarkEvalMPC} are average values and therefore give no direct insight into the real-time feasibility of the MPC solvers and the variation of the computational load over the single sampling steps. To this end, \Fig{\ref{fig:compTimeCumSum}} shows accumulation plots of the computation time per MPC step for three selected problems of the testbench. The computation times were evaluated after 30 successive runs to obtain reliable results. 
The plots show that the computation time of \GRAMPC is almost constant for each MPC iteration, which is important for embedded control applications and to give tight upper bounds on the computation time for real-time guarantees. \ACADO and \VIATOC show a larger variation of the computation time over the iterations, which is mainly due to the active set strategy that both solvers follow and the varying number of QP iterations in each real-time iteration of \ACADO, c.f.~\cite{Houska2011}.

In conclusion, it can be said that \GRAMPC has overall fast and real-time feasible computation times for all benchmark problems, in particular for nonlinear systems and in connection with (nonlinear) constraints. Furthermore, \GRAMPC scales well with increasing system dimension. 

\subsubsection{Embedded realization}\label{sec:embdded}

\begin{table}[t]
	\caption{Computation time and memory usage for the embedded realization of \GRAMPC on \Dspace hardware (DS1202) and ECU level (\Renesas RH850/P1M) with single floating point precision.}
	\label{tab:Embedded}
        \vspace{-2mm}
        	\setlength{\tabcolsep}{3pt}
	\renewcommand*{\arraystretch}{1.1}
	\begin{adjustbox}{max width=\textwidth}
		\begin{tabular}{l r r r | l r r r}
			\toprule
			Problem & $t_\text{\Dspace}$& $t_\text{ECU}$ & Memory & Problem & $t_\text{\Dspace} $  & $t_\text{ECU} $ & Memory \\
			\midrule
			Mass-spring-damper  		& \SI{0.24}{\milli\second} & \SI{4.00}{\milli\second}	& \SI{6.5}{\kilo\byte} & Helicopter    & \SI{0.11}{\milli\second} & \SI{1.88}{\milli\second}	& \SI{5.8}{\kilo\byte} \\
			Motor (PMSM)        		& \SI{0.13}{\milli\second} & \SI{2.10}{\milli\second}	& \SI{2.4}{\kilo\byte} & Quadrotor 	   & \SI{0.21}{\milli\second} & \SI{1.60}{\milli\second}	& \SI{4.4}{\kilo\byte} \\
			Nonl. chain ($m\!=\!4$)		& \SI{4.79}{\milli\second} & $-$ 						& \SI{12.9}{\kilo\byte}& VTOL 		   & \SI{0.37}{\milli\second} & \SI{2.80}{\milli\second}	& \SI{6.2}{\kilo\byte} \\
			Nonl. chain ($m\!=\!6$)		& \SI{9.62}{\milli\second} & $-$ 						& \SI{18.5}{\kilo\byte}& Ball-on-plate & \SI{0.05}{\milli\second} & \SI{0.92}{\milli\second}	& \SI{2.0}{\kilo\byte} \\
			Nonl. chain ($m\!=\!8$)		& \SI{17.50}{\milli\second} & $-$						& \SI{24.1}{\kilo\byte}& Vehicle	   & \SI{0.25}{\milli\second} & \SI{2.69}{\milli\second}	& \SI{3.6}{\kilo\byte} \\
			Nonl. chain ($m\!=\!10$)	& \SI{24.20}{\milli\second} & $-$ 						& \SI{29.8}{\kilo\byte}& CSTR reactor  & \SI{0.43}{\milli\second} & \SI{6.81}{\milli\second}	& \SI{5.1}{\kilo\byte} \\
			2D-Crane 				& \SI{0.18}{\milli\second} & \SI{1.65}{\milli\second}	& \SI{4.5}{\kilo\byte} & PDE reactor   & \SI{6.48}{\milli\second} & $-$ 						& \SI{15.0}{\kilo\byte}\\
			3D-Crane 				& \SI{0.31}{\milli\second} & \SI{3.05}{\milli\second}	& \SI{7.3}{\kilo\byte} &   	  		   &						   & 							&  		   			   \\
			
			\bottomrule 
		\end{tabular}
	\end{adjustbox}
\end{table}

In addition to the general MPC evaluation, this section evaluates the computation time and memory requirements of \GRAMPC for the benchmark problems on real-time and embedded hardware.  \GRAMPC was implemented on a 
\Dspace MicroLabbox (DS1202) with a \SI{2}{\giga\hertz} Freescale QolQ processor as well as on the microntroller RH850/P1M from \Renesas with a CPU frequency of \SI{160}{\mega\hertz}, 
\SI{2}{\mega\byte} program flash and \SI{128}{\kilo\byte} RAM. This processor is typically used in electronic control units (ECU) in the automotive industry. The \GRAMPC implementation on this microcontroller therefore can be seen as a prototypical ECU realization. 
As it is commonly done in embedded systems, \GRAMPC was implemented using single floating point precision on both systems due to the floating point units of the processors.

\Tab{\ref{tab:Embedded}} lists the computation time and RAM memory footprint of \GRAMPC
on both hardware platforms for the testbench problems in Table~\ref{tab:BenchmarkOverview} and \ref{tab:BenchmarkEvalMPC}. The settings of \GRAMPC are the same as in the previous section, except for the floating point precision. Due to 
the compilation size limit of the ECU compiler ($<10$~kB), the nonlinear chain examples as well as the PDE reactor could not be compiled on the ECU.

The computation times on the \Dspace hardware are below the sampling time for all example problems. The same holds for the ECU implementation, except for the 
2D-crane, the PMSM, and the VTOL example. However, as mentioned before, tuning of the algorithm can further reduce the runtime, as most of the multiplier and penalty update parameters are taken as default. Note that there is no constant scaling factor between the computation times on \Dspace and ECU level, which is probably due to the different realization of the math functions by the respective floating point unit / compiler\footnote{For example software or hardware realization of sine or cosine functions.} on the different hardware. 

The required memory is below \SI{9}{\kilo\byte} for all examples, except for the nonlinear chain and the PDE reactor, which is less than \SI{7}{\percent} of the available RAM on the considered ECU. Although the nonlinear chain and the PDE reactor could not be compiled on the ECU as mentioned above,  the memory usage as well as the computation time increase only linearly with the size of the system (using the same \GRAMPC settings). Overall, the computational speed and the small memory footprint demonstrate the applicability of \GRAMPC for embedded systems. 

\subsection{Application examples}
This section discusses some application examples, 
including a more detailed view on two selected problems from the testbench collection (state constrained and semi-implicit problem), a shrinking horizon MPC application, an equality constrained OCP as well as a moving horizon estimation problem.

\subsubsection{Nonlinear constrained model predictive control}\label{sec:NC_MPC}

\begin{figure}[t]
	 \centering
\begin{subfigure}{0.49\textwidth}
	\def\svgwidth{0.9\columnwidth}
	{{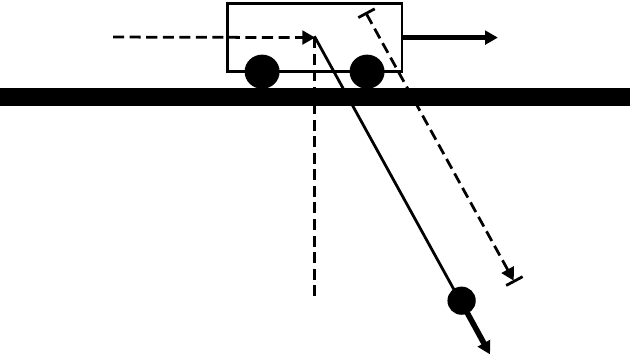}}
\end{subfigure}
\begin{subfigure}{0.49\textwidth}
	\centering
	\vspace{1mm}
\includegraphics[page=1]{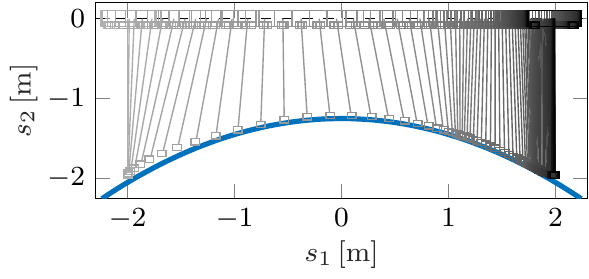}
\end{subfigure}
	\caption{Schematic representation of the overhead crane (left) and 
		simulated crane movement from the initial state to the desired setpoint (right).}
\label{fig:Crane2D}
\end{figure}

The 2D-crane example in Table~\ref{tab:BenchmarkOverview} and \ref{tab:BenchmarkEvalMPC} is a particularly challenging one, as it is subject to a nonlinear constraint that models the collision avoidance of an object or obstacle.
A schematic representation of the overhead crane is given in \Fig{\ref{fig:Crane2D}}. The crane has three degrees-of-freedom and the nonlinear dynamics read as~\cite{Kaepernick2013} 
\begin{align*}
\ddot s_\text{C}
= u_1, \quad 
\ddot s_\text{R}  = u_2, \quad
\ddot \phi = - \frac{1}{s_\text{R}} \left( g \sin(\phi) + a_\text{C} \cos(\phi) + 2 \dot s_\text{R}\dot \phi \right)
\end{align*}
with the gravitational constant $g$. The system state 
$\vm x=[s_\text{C},\dot s_\text{C}, s_\text{R},\dot s_\text{R},\phi,\dot\phi]\TT  \in \mathbb{R}^6$ 
comprises the cart position $s_\text{C}$, the rope length
$s_\text{R}$, the angular deflection $\phi$ and the corresponding
velocities. The two controls $\vm u \in \mathbb{R}^2$ are the cart acceleration $u_1$ and the rope acceleration $u_2$, respectively.

\begin{figure}[t]
	\vspace{1mm}
\includegraphics[page=1]{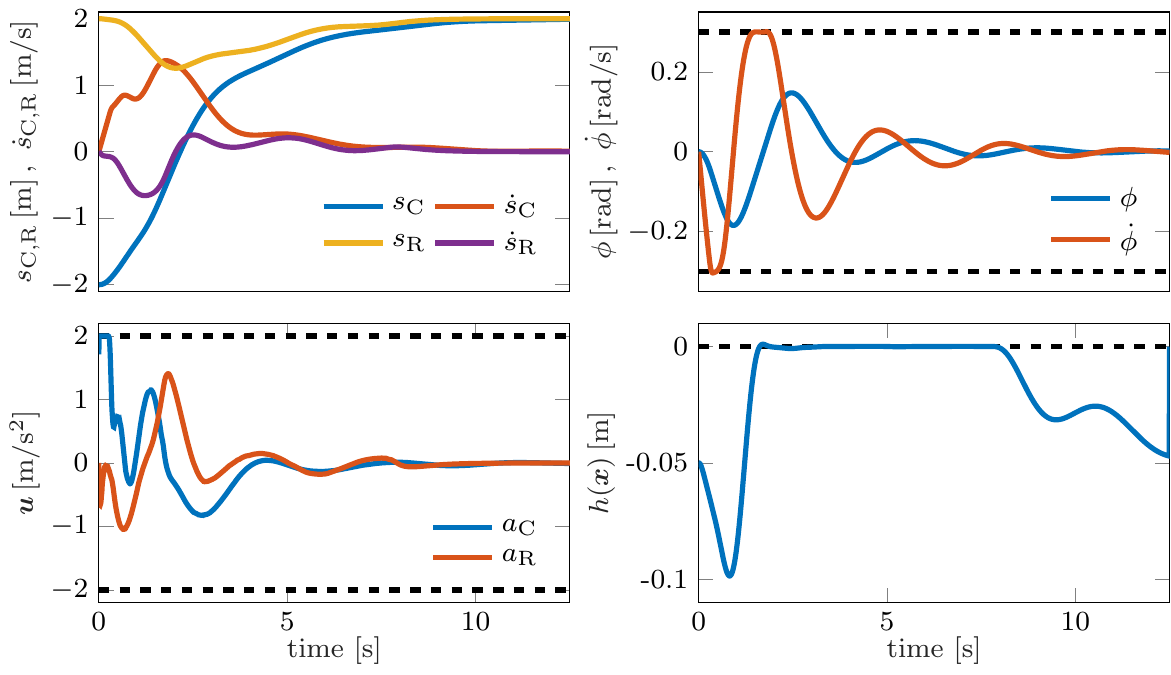}
	\caption{MPC trajectories for the 3DOF crane.}
	\label{fig:CRANE2D_over_obstacle}
\end{figure}

The cost functional~\eqref{eq:MPC_cost} consists of the integral part
\begin{equation}
l(\vm x, \vm u) = (\vm x - \vm x_\text{des} )\TT \vm Q (\vm x - \vm x_\text{des} ) +  
(\vm u - \vm u_\text{des} ) \TT \vm R  (\vm u - \vm u_\text{des} )\,,
\end{equation}
which penalizes the deviation from the desired setpoints $\vm x_\mathrm{des} \in \mathbb{R}^6$ 
and $\vm u_\mathrm{des} \in \mathbb{R}^2$ respectively. The weight matrices are set to 
$\vm{Q}=\diag(1,2,2,1,1,4)$ and $\vm R =\diag(0.05,0.05)$ (omitting units).
The controls and angular velocity are subject to the box constraints
$|u_i| \leq \SI{2}{\meter\per\second\squared},\,i\in{1,2}$ 
and 
$|\dot\phi| \leq \SI{0.3}{\radian\per\second},\,i\in{1,2}$. In addition, the 
nonlinear inequality constraint 
\begin{equation}\label{eq:Crane2D_nlconst}
h (\vm x) = \cos(\phi) s_\text{R} - \SI{0.2}{\per\meter} \left( s_\text{C} + \sin(\phi)  s_\text{R}\right) ^2 + \SI{1.25}{\meter} \leq 0
\end{equation}
is imposed, which represents a geometric security constraint, e.g. for trucks, over which the load 
has to be lifted, see Figure~\ref{fig:CRANE2D_over_obstacle} (right). 

The prediction horizon and sampling time for the crane problem are set to   
$T = \SI{2}{\second}$ and $\Delta t = \SI{2}{\milli\second}$, respectively. The number of augmented Lagrangian steps and inner gradient iterations of \GRAMPC are set to 
$(i_\text{max}, j_\text{max}) = (1,2)$. These settings correspond to the computational results in \Tab{\ref{tab:BenchmarkEvalMPC}} and \ref{tab:Embedded}.

The right part of \Fig{\ref{fig:Crane2D}} illustrates the movement of the overhead crane from the initial 
state $\vm{x}_0 = \left[\SI{-2}{\meter},0, \SI{2}{\meter},0,0,0 \right] \TT$ to the desired setpoint 
$\vm{x}_\text{des} = \left[\SI{2}{\meter},0, \SI{2}{\meter},0,0,0 \right] \TT$. 
\Fig{\ref{fig:CRANE2D_over_obstacle}} shows the corresponding trajectories of the 
states $\vm x(t)$ and controls $\vm u(t)$ as well as the nonlinear 
constraint~\eqref{eq:Crane2D_nlconst} plotted as time function $h(\vm x(t))$.
This transition problem is quite challenging, since the nonlinear constraint \eqref{eq:Crane2D_nlconst} is active for more than half of the simulation time. One can slightly see an initial violation of the constraint $h(\vm x)$ of approximately \SI{1}{\milli\meter}, which should be negligible in practical applications. Nevertheless, one can satisfy the constraint to a higher accuracy by increasing the number of iterations $(i_\text{max}, j_\text{max})$, in particular of the augmented Lagrangian iterations.

\subsubsection{MPC on shrinking horizon}
``Classical'' MPC with a constant horizon length typically 
acts as an asymptotic controller in the sense that a desired setpoint is only reached asymptotically. MPC on a shrinking horizon instead reduces the horizon time $T$ in each sampling step in order to reach the desired setpoint in finite time. In particular, if the desired setpoint is incorporated into a terminal constraint and the prediction horizon $T$ is optimized in each MPC step, then $T$ will be automatically reduced over the runtime due to the principle of optimality.

Shrinking horizon MPC with \GRAMPC is illustrated for the following double integrator problem
\begin{subequations}
\label{eq:shrinkhor:2xint}
\begin{alignat}{3}
\label{eq:shrinkhor:2xint:cost}
& \min_{ u,T} & \quad & J(u,T;\vm x_k) = T +
\frac{1}{2} \int_{0}^{T} r  u^2(\tau) \, \dd \tau 
\hspace{-6cm} &&
\\
& \text{\;s.t.} & \quad & 
\dot{ x}_1(\tau) =  x_2(\tau) \,, \quad 
&&  x_1(0) = x_{1,k} = x_1(t_k)
\\ &&&
\dot{ x}_2(\tau) = u(\tau) \,, \quad
&&  x_2(0) = x_{2,k}=x_2(t_k)
\\
\label{eq:shrinkhor:2xint:boxconstr}
&& \quad & | u(\tau)| \le 1 \,, \quad
&& \tau\in[0,T]
\\ \label{eq:shrinkhor:2xint:terminalconstr}
&& \quad & \vm{ x}(T) = \vm x_\text{des} 
\hspace{-4cm}
\end{alignat}
\end{subequations}
with the state $\vm{ x}=[ x_1, x_2]\TT$ and control $ u$ subject to the box constraint \eqref{eq:shrinkhor:2xint:boxconstr}.
The weight $r>0$ in the cost functional~\eqref{eq:shrinkhor:2xint:cost} 
allows a trade-off
between energy optimality and time optimality of the MPC. The desired setpoint $\vm x_\text{des}$ is added as terminal constraint \eqref{eq:shrinkhor:2xint:terminalconstr}, i.e.\ $\vm g_T(\vm{ x}(T)) := \vm{ x}(T) - \vm x_\text{des} = \mb 0$ in view of~\eqref{eq:OCP_eqconstr}, and the prediction horizon $T$ is treated as optimization variable in addition to the control $ u(\tau)$, $\tau\in[0,T]$.

\begin{figure}[t]
	\includegraphics[page=1]{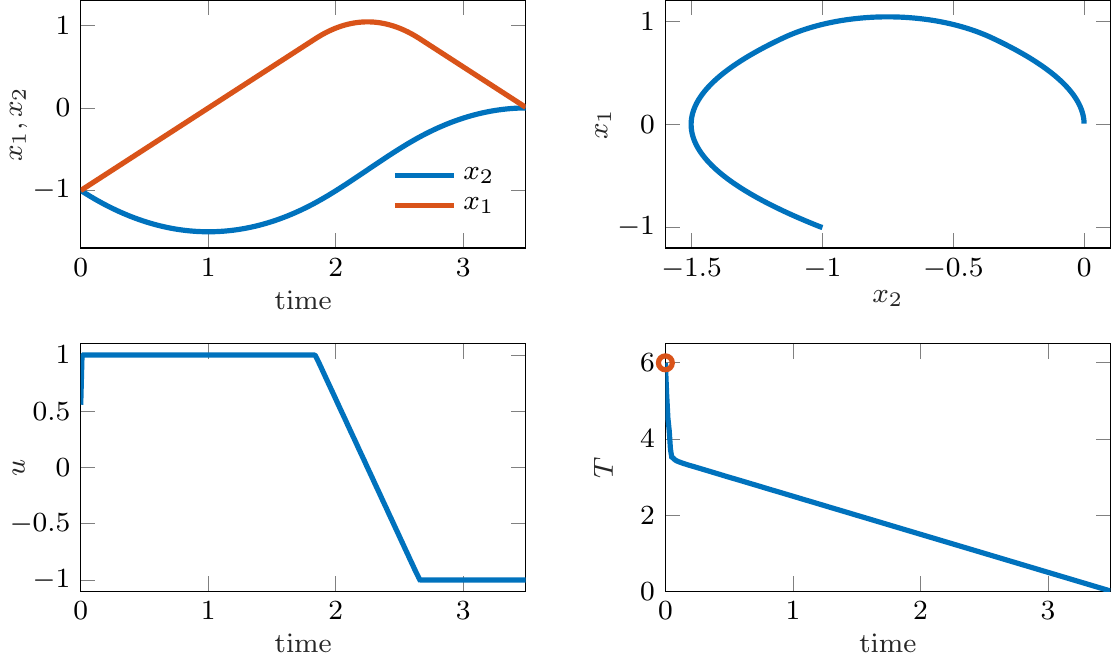}
	\caption{MPC trajectories with shrinking horizon for the double integrator problem~\eqref{eq:shrinkhor:2xint}.}
	\label{fig:DoubleIntegratorTrajectories}
\end{figure}

For the simulations, the weight in the cost functional is set to $r = 0.01$ and the initial value of the horizon length is chosen as $T=6$. The iteration numbers for \GRAMPC are set to $(i_\text{max},j_\text{max})=(1,2)$ in conjunction with a sampling time of $\Delta t=0.001$ in order to resemble a real-time implementation.
\Fig{\ref{fig:DoubleIntegratorTrajectories}} shows the simulation results for the double integrator problem with the desired setpoint $\vm x_\text{des}=\vm 0$ and the initial state $\vm x(0)=[-1,-1]\TT$.
Obviously, the state trajectories reach the origin in finite time corresponding to the anticipated behavior of the shrinking horizon MPC scheme. 
The lower right plot of \Fig{\ref{fig:DoubleIntegratorTrajectories}} shows the temporal evolution of the horizon length $T$ over the runtime $t$. The initial end time of $T=6$ is marked as a red circle. In the first MPC steps, the optimization quickly decreases the end time to approximately $T=3.5$. In this short time interval, \GRAMPC produces a suboptimal solution due to the strict limitation of the iterations $(i_\text{max},j_\text{max})$. Afterwards, however, the prediction horizon declines linearly, which concurs with the principle of optimality and shows the optimality of the computed trajectories after the initialization phase.
In the simulated example, this knowledge is incorporated in the MPC implementation by substracting the sampling time $\Delta t$ from the last solution of $T$ for warm-starting the next MPC step.
The simulation is stopped as soon as the horizon length reaches its minimum value $T_\mathrm{min} = \Delta t = 0.1$.
\subsubsection{Semi-implicit problems}
\label{sec:reactor_PDE}

The system formulations that can be handled with \GRAMPC include DAE systems with singular mass matrix $\vm M$ as well as general semi-implicit systems. 
An application example is the discretization of spatially distributed systems by means of finite elements.
This is illustrated for a quasi-linear diffusion-convection-reaction system, which is also implemented in the testbench (PDE reactor example). 
The thermal behavior of the reactor is described on the one-dimensional spatial domain $z=(0,1)$ using the PDE formulation~\cite{Utz2010}
\begin{subequations}
\label{eq:chap5:SysDynReactor__PDE}
\begin{equation}
\label{eq:chap5:SysDynReactor__PDE_heateq}
\partial_t \theta = \partial_z \left[ (q_0+q_1 \theta) \partial_z \theta -\nu \theta \right] + (r_0+r_1 \theta) \theta	
\end{equation}
with the boundary and initial conditions
\begin{alignat}{2}
\label{eq:chap5:SysDynReactor__PDE_NeumannRB}
\partial_z \theta|_{z=0}  &= 0
\\ \label{eq:chap5:SysDynReactor__PDE_BoundaryCtrl}
\partial_z \theta|_{z=1} + \theta(1,t) &= u
\\
\theta(\cdot,0)&=\theta_0
\end{alignat}
\end{subequations}
for the temperature $\theta=\theta(z,t)$. The process is controlled by the boundary control $u=u(t)$. 
Diffusive and convective processes of the reactor are modeled by the nonlinear heat equation \eqref{eq:chap5:SysDynReactor__PDE_heateq} with the parameters 
$q_1=\num{2}$, $q_2=\num{-0.05}$, and $\nu=\num{1}$, respectively. 
Reaction effects are included using the parameters $r_0=\num{1}$ and $r_1=\num{0.2}$. 
The Neumann boundary condition \eqref{eq:chap5:SysDynReactor__PDE_NeumannRB}, the Robin boundary condition \eqref{eq:chap5:SysDynReactor__PDE_BoundaryCtrl}, and the initial condition \eqref{eq:chap5:SysDynReactor__PDE} complete the mathematical description of the system dynamics. 
Both spatial and time domain are normalized for the sake of simplicity.
A more detailed description of the system dynamics can be found in~\cite{Utz2010}.

\begin{figure}[tb!]
	\includegraphics[page=1]{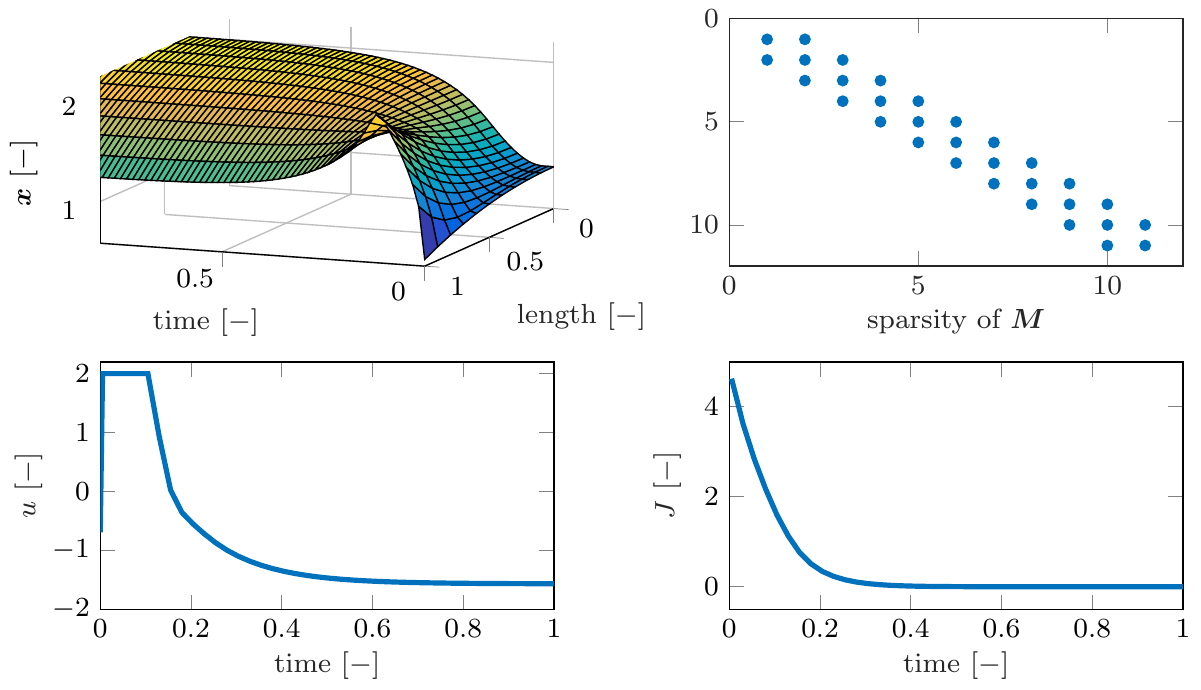}
	\caption{Simulated MPC trajectories for the PDE reactor~\eqref{eq:chap5:SysDynReactor__PDE}.}
	\label{fig:PDEReactor}
\end{figure}

The PDE system \eqref{eq:chap5:SysDynReactor__PDE} is approximated by an ODE system of the form \eqref{eq:optxcond} by applying a finite element discretization technique~\cite{zienkiewicz_book1983_fem}, whereby the new state variables $\vm x\in\mathbb{R}^{N_{\vm x}}$ approximate the temperature $\theta$ on the discretized spatial domain $z$ with $N_{\vm x}$ spatial grid points.
The finite element discretization eventually leads to a finite-dimensional system dynamics of the form
$\vm M\vm{\dot x} = \vm f(\vm x,u)$
with the mass matrix $\vm M\in\mathbb{R}^{N_{\vm x}\times N_{\vm x}}$ and the nonlinear system function $\vm f(\vm x, u)$.
In particular, $N_{\vm x}=11$ grid points are chosen for the \GRAMPC simulation of the  reactor~\eqref{eq:chap5:SysDynReactor__PDE}.
The upper right plot of Figure~\ref{fig:PDEReactor} shows the sparsity structure 
of the mass matrix $\vm M\in\mathbb R^{11\times 11}$.

The control task for the reactor is the stabilization of a stationary profile $\vm x_\text{des}$ 
which is accounted for in the quadratic MPC cost functional
\begin{align*}
J(u;\vm x_k) := \tfrac{1}{2}\big\|\vm x(T) - \vm x_{\rm des}\big\|^2
+ \int_0^T \tfrac{1}{2}\big\|\vm x(t) - \vm x_{\rm des}\big\|^2 + \num{e-2} \big\|u-u_{\rm des}\big\|^2
\,\dd t\,.
\end{align*}
The desired setpoint $(\vm x_\text{des},u_\text{des})$ as well as the initial values $(\vm x_0,u_0)$ are numerically determined from the stationary differential equation
\begin{alignat}{2}
0 &= \partial_z \left[ (q_0+q_1
\theta(z,\tau)) \partial_z \theta(z,\tau) -\nu \theta(z,\tau) \right] 
+ (r_0+r_1 \theta(z,\tau)) \theta(z,\tau)
\end{alignat}
with the corresponding boundary conditions
\begin{alignat*}{2}
\theta(0,0) & = 1, &\quad \partial_z \theta(0,0) & = 0 \\
\theta(0,\infty) & = 2, &\quad \partial_z \theta(0,\infty) & = 0\,.
\end{alignat*}
The prediction horizon and sampling time of the MPC scheme are set to $T=\num{0.2}$ 
and $\Delta t=\num{0.005}$. The number of iterations are limited by $(i_\text{max},j_\text{max})=(1,2)$. The box constraints for the control are chosen as $|u(t)| \le 2$. 

The numerical integration in \GRAMPC is carried out using the solver RODAS~\cite{Hairer:Book:1996:Stiff}. The sparse numerics of RODAS allow one to cope with the banded structure of the matrices in a computationally efficient manner. \Fig{\ref{fig:PDEReactor}} shows the setpoint transition from the initial temperature profile $\vm x_0$ to the desired temperature $\vm x_{\rm des}=[2.00,    1.99,    1.97,    1.93,    1.88,    1.81,    1.73,    1.63,    1.51,    1.38,    1.23]\TT$ and desired control $u_\text{des}=-1.57$. 

\subsubsection{OCP with equality constraints}
An illustrative example of an optimal control problem with equality constraints is a dual arm robot with closed kinematics, e.g.\ to handle or carry work pieces with both arms. 
For simplicity, a planar dual arm robot with the joint angles $(x_1, x_2, x_3)$ and $(x_4, x_5, x_6)$ for the left and right arm is considered. 
The dynamics are given by simple integrators
\begin{align}
\vm{\dot x} = \begin{bmatrix}
\dot x_1 & \dot x_2 & \dot x_3 & \dot x_4 & \dot x_5 & \dot x_6
\end{bmatrix}\TT
= \begin{bmatrix}
u_1 & u_2 & u_3 & u_4 & u_5 & u_6
\end{bmatrix}\TT
= \vm u
\end{align}
with the joint velocities as control input $\vm u$. 
Given the link lengths $\vm a = [a_1$, $a_2$, $a_3$, $a_4$, $a_5$, $a_6]$, the forward kinematics of left and right arm are computed by
\begin{align}
\vm p_{\mathrm L}(\vm x) = \begin{bmatrix}
a_1 \cos(x_1) + a_2 \cos(x_1 + x_2) + a_3 \cos(x_1 + x_2 + x_3) \\
a_1 \sin(x_1) + a_2 \sin(x_1 + x_2) + a_3 \sin(x_1 + x_2 + x_3) \\
x_1 + x_2 + x_3
\end{bmatrix}
\end{align}
and
\begin{align}
\vm p_{\mathrm R}(\vm x) = \begin{bmatrix}
1 + a_4 \cos(x_4) + a_5 \cos(x_4 + x_5) + a_6 \cos(x_4 + x_5 + x_6) \\
a_4 \sin(x_4) + a_5 \sin(x_4 + x_5) + a_6 \sin(x_4 + x_5 + x_6) \\
x_4 + x_5 + x_6
\end{bmatrix} \,.
\end{align} 
The closed kinematic chain is enforced by the equality constraint
\begin{align}
\vm g(\vm x) := \vm p_{\mathrm L}(\vm x) - \vm p_{\mathrm R}(\vm x) - \left[0,\, 0,\, \pi \right]\TT = \vm 0 \,.
\end{align}
A point-to-point motion from $\vm x_0 = [ \frac{\pi}{2}, -\frac{\pi}{2}, 0, -\frac{\pi}{2}, \frac{\pi}{2}, 0 ]\TT$ to 
$\vm x_f = [ -\frac{\pi}{2}$, $\frac{\pi}{2}$, $0$, $\frac{\pi}{2}$, $-\frac{\pi}{2} ]\TT$ is considered as control task, which is formulated as the optimal control problem
\begin{subequations}
\begin{alignat}{3}
&& \min_{\vm u} \quad& J(\vm u) := \int_0^T \frac12 \vm u(t)\TT \vm R \vm u(t) \, \dd t 
\hspace{-4cm}
\\
&& \text{s.t.} \quad& \vm{\dot x}(t) = \vm u(t) \,,\quad 
&& \vm x(0) = \vm x_0 \,,\quad \vm x(T) = \vm x_f \\
&& & \vm g(\vm x(t)) = \vm 0 \,,\quad && \vm u(t) \in \left[ \vm u_\text{min}, \vm u_\text{max} \right]
\end{alignat}\label{eq:BenchmarkEvalRobotOCP}
\end{subequations}
with the fixed end time $T = 10 \, \text{s}$ and the control constraints 
$-\vm u_\text{min} = \vm u_\text{max} = [1,1,1,1,1,1]\TT\text{s}^{-1}$, 
which limit the angular speeds of the robot arms. The cost functional minimizes the squared joint velocities with $\vm R = \eye{6}$.

\begin{table}[t]
	\caption{Computation results for the planar two-arm robot with closed kinematics.}
	\label{tab:BenchmarkEvalRobotOCP}
        \vspace{-2mm}
	\centering
	\begin{adjustbox}{max width=\textwidth}
		\begin{tabular}{@{}p{0.2cm}cccccp{0cm}}
		\toprule
		& Gradient & Constraint & Gradient iter- & Augm.\ Lagr. & $t_\text{CPU}$ & \\
                & tol.~$\varepsilon_\text{rel,c}$ & tol.~$\vm\varepsilon_{\vm g}$ & ations $i$ (avg.) & iterations $j$ & [ms] \\
		\midrule
		& $\num{1e-5}$ & $\num{1e-3}$ & 64 & 189 & $135$ & \\
		& $\num{1e-6}$ & $\num{1e-4}$ & 65 & 298 & $214$ & \\
		& $\num{1e-7}$ & $\num{1e-5}$ & 306 & 190 & $628$ & \\
		& $\num{1e-8}$ & $\num{1e-6}$ & 222 & 431 & $1015$ & \\
		\bottomrule
		\end{tabular}
	\end{adjustbox}
\end{table}

Table~\ref{tab:BenchmarkEvalRobotOCP} shows the computation results of \GRAMPC for solving OCP~\eqref{eq:BenchmarkEvalRobotOCP} with increasingly restrictive values of the gradient tolerance $\varepsilon_\text{rel,c}$ and constraint tolerance $\vm\varepsilon_{\vm g}$ that are used for checking convergence of Algorithm~1 and 2.
The required computation time $t_\text{CPU}$ as well as the average number of (inner) gradient iterations and the number of (outer) augmented Lagrangian iterations are shown in \Tab{\ref{tab:BenchmarkEvalRobotOCP}}. 
The successive reduction of the tolerances $\varepsilon_\text{rel,c}$ and $\vm\varepsilon_{\vm g}$
highlights that the augmented Lagrangian framework is able to quickly compute a solution with moderate accuracy. When further tightening the tolerances, the computation time as well as the required iterations increase clearly. This is to be expected as augmented Lagrangian and gradient algorithms are linearly convergent opposed to super-linear or quadratic convergence of, e.g., SQP or interior point methods. However, as the main application of \GRAMPC is model predictive control, the ability to quickly compute suboptimal solutions that are improved over time is more important than the numerical solution for very small tolerances.

\begin{figure}[t]
	\includegraphics[page=1]{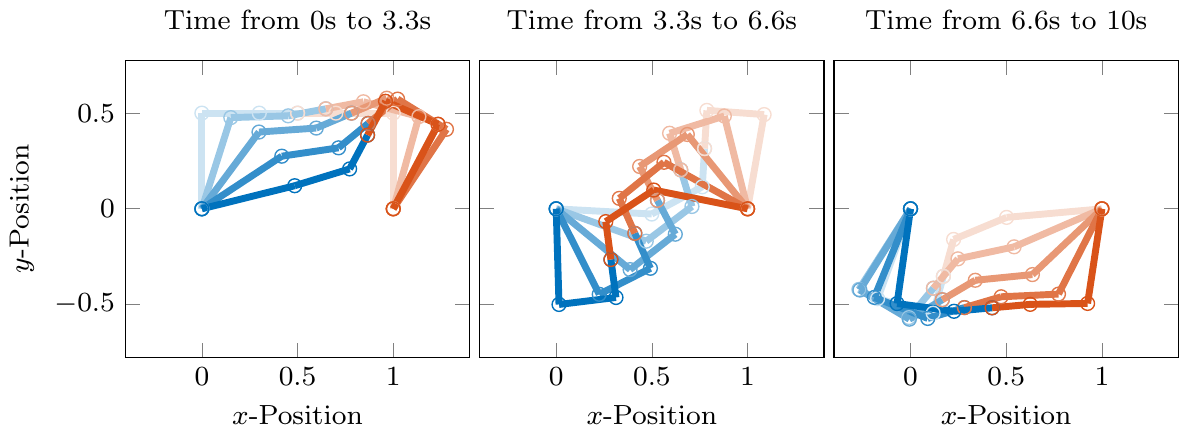}
	\caption{Solution trajectory for the planar two-arm robot with closed kinematics.}
	\label{fig:BenchmarkEvalRobotOCP}
\end{figure}

The resulting trajectory of the planar robot is depicted in \Fig{\ref{fig:BenchmarkEvalRobotOCP}} and shows that the solution of the optimal control problem~\eqref{eq:BenchmarkEvalRobotOCP} involves moving through singular configurations of left and right arm, which makes this problem quite challenging.

\subsubsection{Moving horizon estimation}

Another application domain for \GRAMPC is moving horizon estimation (MHE) by taking advantage of the parameter optimization functionality. This is illustrated for the CSTR reactor model listed in the MPC testbench in Table~\ref{tab:BenchmarkOverview}.
The system dynamics of the reactor is given by~\cite{Rothfuss1996}
\begin{subequations}
\begin{align}
  \dot c_\mathrm{A} &= -k_1(T) c_\mathrm{A} - k_2(T) c_\mathrm{A}^2 + (c_{\mathrm{in}} - c_\mathrm{A}) u_1 \\
  \dot c_\mathrm{B} &= k_1(T) c_\mathrm{A} - k_1(T) c_\mathrm{B} - c_\mathrm{B} u_1 \\
  \dot{T}~ &= -\delta ( k_1(T) c_\mathrm{A} \Delta H_\mathrm{AB} + k_1(T) \Delta H_\mathrm{BC} + k_2(T) c_\mathrm{A}^2 \Delta H_\mathrm{AD}) \nonumber\\
                       &\hspace{1cm}+ \alpha (T_\mathrm{C} - T) + (T_\mathrm{in} - T) u_1 \\
  \dot{T}_\mathrm{C} &= \beta (T - T_\mathrm{C}) + \gamma u_2
\end{align}
\end{subequations}
with the state vector $\vm x=[c_\mathrm{A}, c_\mathrm{B},T, T_\mathrm{C}]\TT$ consisting of the monomer and product concentrations $c_\mathrm{A}$, $c_\mathrm{B}$ and the reactor and cooling temperature $T$ and $T_\mathrm{C}$. The control variables $\vm u=[u_1,u_2]\TT$ are the normalized flow rate and cooling power. The measured quantities are the temperatures
$y_1 = T$ and $y_2 = T_\mathrm{C}$.
The parameter values and a more detailed description of the system are given in \cite{Rothfuss1996}.

\begin{figure}[t]
	\includegraphics[page=1]{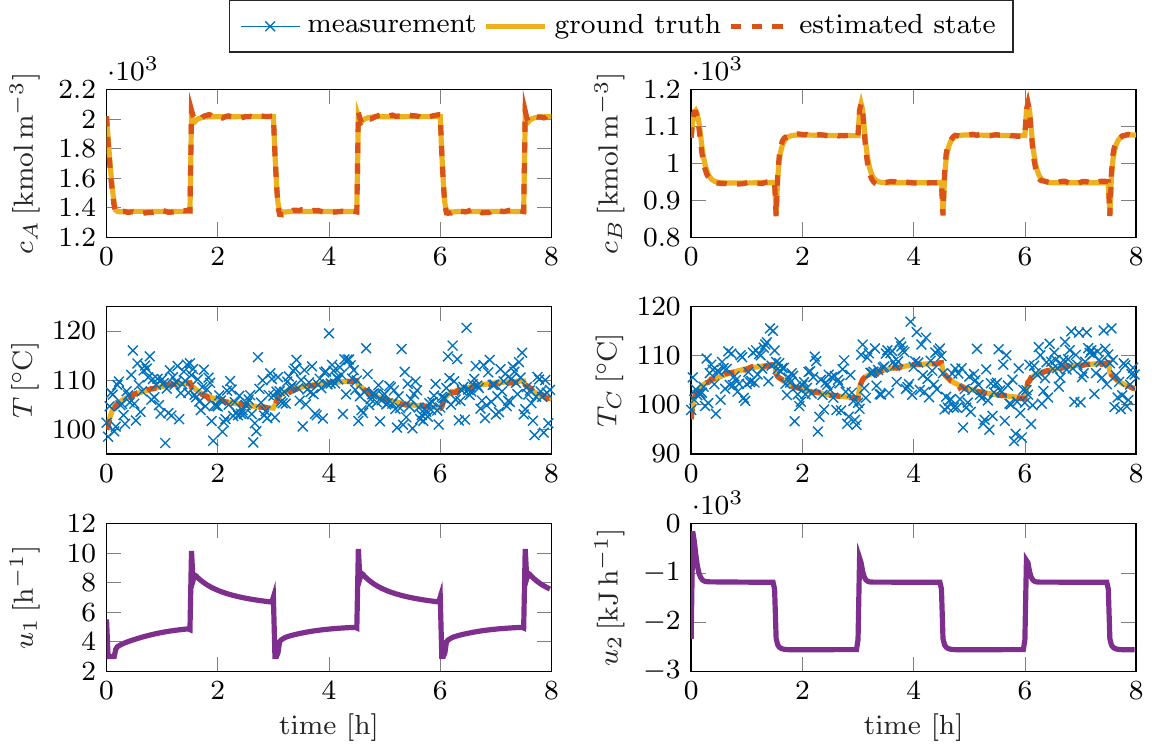}
	\caption{MHE/MPC simulation results for the CSTR reactor.}
	\label{fig:MHEexample}
\end{figure}

The following scenario considers the interconnection of the MHE with the MPC from the testbench, i.e.\
the state $\vm{\hat x}_k$ at the sampling time $t_k$ is estimated and provided to the MPC.
The cost functional of the MPC is designed according to
\begin{align}
J(\vm u;\vm x_k) := \Delta \mb x^\mathrm{T}(T) \mb P \Delta \mb x(T)
+ \int_0^T \Delta \mb x^\mathrm{T}(t) \mb Q \Delta \mb x(t) + \Delta \mb u^\mathrm{T}(t) \mb R \Delta \mb u(t)\,\dd t\,,
\end{align}
where $\Delta \mb x = \mb x - \mb x_{\mathrm{des}}$ and $\Delta \mb u = \mb u - \mb u_{\mathrm{des}}$ penalize the deviation of the state and control from the desired setpoint $(\mb x_{\mathrm{des}},\mb u_{\mathrm{des}})$.
The MHE uses the cost functional defined in Section~\ref{sec:MHE}, c.f. \eqref{eq:MHE_1_cost}. The sampling rate for both MPC and MHE is set to
{$\Delta t=\SI{1}{\second}$}.
The MPC runs with a prediction horizon of 
{$T=\SI{20}{\minute}$}
and $40$ discretization points, while the MHE horizon is set to 
{$T_\mathrm{MHE}=\SI{10}{\second}$} with $10$ discretization points. The \GRAMPC implementation of the MHE uses only one gradient iteration per sampling step, i.e.\ $(i_\text{max},j_\text{max})=(1,1)$, while the implementation of the MPC uses three gradient iterations, i.e.\ $(i_\text{max},j_\text{max}) =(1,3)$.

The simulated scenario in Figure~\ref{fig:MHEexample} consists of alternating setpoint changes between two stationary points. The two measured temperatures are subject to a Gaussian measurement noise with a standard deviation of $\SI{4}{\degreeCelsius}$. The initial guess for the state vector $\vm x$ of the MHE differs from the true initial values by $\SI{100}{\kilo\mol\per\cubic\meter}$ for both concentrations $c_\text{A}$ and $c_\text{B}$ and by $\SI{5}{\degreeCelsius}$, respectively $\SI{-7}{\degreeCelsius}$, for the reactor and cooling temperature. This initial disturbance vanishes within a few iterations and the MHE
tracks the ground truth (i.e. the simulated state values) closely, with an average error of $\delta\vm{\hat x} =[\SI{7.12}{\kilo\mol\per\cubic\meter}$, $\SI{6.24}{\kilo\mol\per\cubic\meter}$, $\SI{0.10}{\degreeCelsius}$, $\SI{0.09}{\degreeCelsius}]\TT$.
This corresponds to a relative error of less than $\SI{0.1}{\percent}$ for each state variable, when normalized to the respective maximum value.
One iteration of the MHE requires a computation time of $t_\mathrm{CPU}=\SI{11}{\micro\second}$ to calculate a new estimate of the
state vector $\vm x$ and therefore about one third of the time necessary for one MPC iteration.

\section{Conclusions}
\label{sec:conclusions}

The augmented Lagrangian algorithm presented in this paper
is implemented in the nonlinear model predictive control (MPC) toolbox
\GRAMPC and extends its original version that was published in 2014 in several
significant aspects. The system class that can be handled by \GRAMPC
are general nonlinear systems described by 
explicit or semi-implicit differential equations or
differential-algebraic equations (DAE) of index~1. Besides input
constraints, the algorithm accounts 
for nonlinear state and/or input dependent equality and inequality
constraints as well as for unknown parameters and a possibly free end
time as further optimization variables, which is relevant, for instance, for moving
horizon estimation or MPC on a shrinking horizon. 
The computational efficiency of \GRAMPC~is
illustrated for a testbench of representative MPC problems and in
comparison with two state-of-the-art nonlinear MPC toolboxes. In
particular, the augmented Lagrangian algorithm implemented in \GRAMPC
is tailored to embedded MPC with very low memory requirements. This is
demonstrated in terms of runtime results on dSPACE and ECU hardware
that is typically used in automotive applications. 
\GRAMPC is available at \url{http://sourceforge.net/projects/grampc}
and is licensed under the GNU Lesser General Public License
(version~3), which is suitable for both academic and
industrial/commercial purposes. 
\begin{appendix}
\section{Transformation of inequality to equality constraints}\label{sec:AppIeqc}
This appendix describes the transformation of the inequality constraints \eqref{eq:OCP_ieqconstr} 
to the equality constraints~\eqref{eq:transfieqc} in more detail.
A standard approach of augmented Lagrangian methods is to introduce 
slack variables $ \vm v \geq 0$ and $ \vm v_T \geq 0 $ in order to
add \eqref{eq:OCP_ieqconstr} to the existing equality constraints
\eqref{eq:OCP_eqconstr} according to 
\begin{equation}
\label{eq:eqconstr_stacked_1}
\vm{\hat g}_T(\vm x, \vm p, T, \vm v_T) = \begin{bmatrix}
\vm g_T(\vm x, \vm p, T) \\
\vm h_T(\vm x, \vm p, T) + \vm v_T
\end{bmatrix} = \vm 0
\,,
\quad
\vm{\hat g}(\vm x, \vm u, \vm p, t, \vm v) = \begin{bmatrix}
\vm g(\vm x, \vm u, \vm p, t) \\
\vm h(\vm x, \vm u, \vm p, t) + \vm v
\end{bmatrix} = \vm 0\,.
\end{equation}
The set of constraints~\eqref{eq:eqconstr_stacked_1} are then
adjoined to the cost functional~\eqref{eq:OCP_cost}
\begin{equation}
\label{eq:Appcost_augm_1}
\hat J(\vm u, \vm p, T, \vm \lag, \vm \lag_T, \vm\pen,\vm \pen_T, \vm v, \vm v_T;\vm x_0) =
\hat V(\vm x, \vm p, T, \vm \lag_T, \vm \pen_T, \vm v_T) 
+ \int_0^T \hat l(\vm x, \vm u, \vm p, t, \vm \lag, \vm \pen, \vm v) \, \mathrm dt
\end{equation}
with the new terminal and integral cost functions
\begin{subequations}
	\label{eq:costII_V&L}
	\begin{align}
	\hat V(\vm x, \vm p, T, \vm \lag_T, \vm \pen_T, \vm v_T) 
	= &V(\vm x, \vm p, T) + \vm\mu_T\TT\vm{\hat g}_T(\vm x, \vm p, T, \vm v_T) 
	+ \frac12 \big\| \vm{\hat g}_T(\vm x, \vm p, T,\vm v_T) \big\|^2_{\vm \penmatrix_T}
	\\
	\hat l(\vm x, \vm u, \vm p, t, \vm \lag, \vm \pen, \vm v) 
	= &l(\vm x, \vm u, \vm p, t)
	+ \vm\mu\TT\vm{\hat g}(\vm x, \vm p, \vm u,t, \vm v_T)
	+ \frac12 \big\| \vm{\hat g}(\vm x, \vm p, \vm u,t,\vm v_T) \big\|^2_{\vm \penmatrix}
	\end{align}
\end{subequations}
using the multipliers $ \vm \lag $ and $ \vm \lag_T $, 
the penalties $ \vm \pen  $ and $ \vm \pen_T$ and 
the corresponding diagonal matrices $\vm \penmatrix$ and $\vm \penmatrix_T$.
Instead of solving the original OCP~\eqref{eq:OCP_orig}, the augmented
Lagrangian approach solves the max-min-problem
\begin{subequations}\label{eq:OCPII}
	\begin{align}
	\max_{\vm \lag, \vm \lag_T} \, \min_{\vm u, \vm p, T, \vm v, \vm v_T} \quad& 
	\hat J(\vm u, \vm p, T, \vm \lag, \vm \lag_T, \vm \pen, \vm \pen_T,\vm v,\vm v_T;\vm x_0) 
	%= 
	%\hat V(\vm x, \vm p, T, \vm \lag_T, \vm \pen_T, \vm v_T) 
	%+ \int_0\TT \hat l(\vm x, \vm u, \vm p, t, \vm \lag, \vm c, \vm v) \, \mathrm dt
	\label{eq:OCP2_cost}
	\\ \label{eq:OCP2_dyn}
	\textrm{s.t.~~~~} \quad& \vm M \vm{\dot x}(t) = \vm f(\vm x, \vm u, \vm p, t) 
	\,,\quad 
	\vm x(0) = \vm x_0%\,,\quad  t \in [0, T]
	\\
	& \vm v(t) \geq 0\,, \quad \vm v_T \geq 0\,,\quad \vm u(t) \in [\vm u_{\min}, \vm u_{\max}] \\%\,,\quad  t \in [0, T]\\
	& 
	\vm p \in [\vm p_{\min}, \vm p_{\max}]
	\,,\quad 
	T \in [T_{\min}, T_{\max}] \,.
	\end{align}
\end{subequations}
The minimization of the cost functional~\eqref{eq:OCP2_cost} with
respect to $\vm v$ and $\vm v_T$ can be carried out explicitly. %~\cite{Bertsekas1996}.  
In case of $\vm v_T$, the minimization of \eqref{eq:OCP2_cost} reduces
to the strictly convex, quadratic problem
\begin{equation}
\label{eq:min_slack_T}
\min_{\vm v_T\ge \vm 0} ~~ \lag_{\vm{h}_T}\TT \big( \vm h_T(\vm x, \vm p, T) + \vm v_T\big) + \frac12 \| \vm h_T(\vm x, \vm p, T) + \vm v_T \|^2_{\vm \penmatrix_{\vm{h}_T}}
\end{equation}
for which the solution follows from the stationarity condition and projection if
the constraint $\vm v_T\ge \vm 0$ is violated
\begin{equation}
\label{eq:slack_opt_T}
\vm v_T = \vm\max\left\{ \vm 0,-\vm\penmatrix_{\vm h_T}^{-1}\vm \lag_{\vm h_T}- \vm h_T(\vm x, \vm p, T)\right\}\,.
\end{equation}
The minimization of \eqref{eq:OCP2_cost} w.r.t.\ the slack
variables $\vm v=\vm v(t)$ corresponds to
\begin{equation}
\label{eq:min_slack}
\min_{\vm v(\cdot)\ge \vm 0} ~~ \int_0^T \vm\lag_{\vm{h}}\TT \big( \vm h(\vm x,
\vm u, \vm p, t) + \vm v \big)  + \frac12 \| \vm h(\vm x, \vm u, \vm
p, t) + \vm v \|^2_{\vm \penmatrix_{\vm{h}}} \dd t\,.
\end{equation}
Since $\vm v$ only occurs in the integral and is not influenced by 
the dynamics~\eqref{eq:OCP2_dyn}, the minimization of \eqref{eq:OCP2_cost}
w.r.t.\ $\vm v=\vm v(t)$ can be carried out pointwise in time and
therefore reduces to a convex, quadratic problem similar to \eqref{eq:min_slack_T}
with the pointwise solution
\begin{equation}
\label{eq:slack_opt}
\vm v = \vm\max\left\{\vm 0,-\vm\penmatrix_{\vm h}^{-1}\vm \lag_{\vm h}-\vm h(\vm x, \vm u, \vm p, t)\right\}\,.
\end{equation}
Inserting the solutions \eqref{eq:slack_opt_T} and
\eqref{eq:slack_opt} into~\eqref{eq:eqconstr_stacked_1} eventually yields the 
transformed equality constraints~\eqref{eq:transfieqc}.

\end{appendix}
\bibliographystyle{spmpsci}
\bibliography{grampc}   
\end{document}